\documentclass[11pt]{article}
\usepackage{amssymb, amsmath, amsthm}
\usepackage{graphicx}
\usepackage{color}
\newtheorem{thm}{Theorem}[section]
\newtheorem{lem}{Lemma}[section]
\newtheorem{cor}[lem]{Corollary}

\theoremstyle{definition}
\newtheorem{defn}[thm]{Definition}
\theoremstyle{remark}
\newtheorem{rem}[thm]{Remark}
\numberwithin{equation}{section}

\newcommand{\E}{\mathbf{E}\,}

\newcommand{\Tr}{\mathrm{Tr}\;\!}

\newcommand{\re}{\mathrm{Re}\;\!}
\newcommand{\im}{\mathrm{Im}\;\!}
\newcommand{\mc}{}
\newenvironment{Proof of}{\removelastskip\par\medskip
\noindent{\em Proof of} \rm}{\penalty-20\null\hfill$\square$\par\medbreak}

\begin{document}
\setcounter{page}{1}

\title{\bf On the Circular Law}

\author{{\bf F. G\"otze}$^{1}$\\{\small Faculty of Mathematics}
\\{\small University of Bielefeld}\\{\small Germany}
\and {\bf A. Tikhomirov}$^{1,2}$\\{\small Faculty of Mathematics}
\\{\small University of Bielefeld}\\{\small Germany}
\\{\small Faculty of Mathematics and Mechanics}\\{\small Sankt-Peterburg State University}
\\{\small S.-Peterbudg, Russia}}
\date{}
\maketitle
 \footnote{$^1$Partially supported by INTAS grant  N
03-51-5018.} \footnote{$^2$Partially supported by RFBF--DFG
grant N
04-01-04000, by RF grant of the leading scientific schools NSh-4222.2006.1}

\today

\begin{abstract}
We consider the joint distribution of real and imaginary parts of
eigenvalues of random matrices with independent real entries with mean
zero and unit variance. We prove the convergence of this distribution to
the uniform distribution on the unit  disc without assumptions on the
existence of a  density for the distribution of entries. We assume however
that the entries have sub-Gaussian tails or are sparsely non-zero.
\end{abstract}

\maketitle
\markboth{ F. G\"otze, A.Tikhomirov}{On the circular law}

\section{Introduction}

Let $X_{jk},{  }1\le j, k<\infty$,
be complex random variables
with $\mathbb E X_{jk}=0$ and
$\mathbb E |X_{jk}|^2=1$.
For a fixed $n\geq 1$,
 denote by $\lambda_1,\ldots,\lambda_n$ the eigenvalues of
the  $n\times n$ matrix
\begin{equation}\label{sym}
\bold X=\frac1{\sqrt n}(X_n{(j,k)})_{j,k=1}^n,
\quad X_n{(j,k)}=\frac1{\sqrt n}X_{jk},\text{  for  }1\leq j,k\leq n,
\end{equation}
and define its empirical  spectral distribution function by
\begin{equation}
G_n(x,y)=\frac1n\sum_{j=1}^nI_{\{\re\{\lambda_j\}\leq x,\;\im\{\lambda_j\}\le y\}},
\end{equation}
where $I_{\{B\}}$ denotes the indicator of an event $B$. We investigate
 the convergence of the expected spectral distribution function
$\mathbb E G_n(x,y)$ to the distribution function $G(x,y)$ of the uniform
distribution over the unit disc in $\mathbb R^2$.

We shall assume that the random variables $X_{jk}$ are sub-Gaussian, i. e.
\begin{defn} A random variable $\beta$ is called sub-Gaussian (respectively 
$\beta$ has a distribution with
 sub-Gaussian tails)  if for any $t>0$
\begin{equation}\notag
\Pr\{|\beta|>t\}\le C\exp\{-ct^2\}.
\end{equation}
\end{defn}

The main result of our paper is the following

\begin{thm}\label{thm0}
Let $X_{jk}$ be independent identically distributed  sub-Gaussian random variables with
\begin{equation}\notag
\E X_{jk}=0,\qquad \E |X_{jk}|^2=1.
\end{equation}
Then  $\E G_n(x,y)$ converges weakly to the distribution function $G(x,y)$ as $n\to\infty$.
\end{thm}

We shall prove the same result for the follows class of sparse matrices. Let
$\varepsilon_{jk}$, $j,k=1,\ldots,n$  denote Bernoulli random variables which are independent in aggregate and
independent of $(X_{jk})_{j,k=1}^n$  with
$p_n:=\Pr\{\varepsilon_{jk}=1\}$. Consider the matrix $\bold
X^{(\varepsilon)}=\frac1{\sqrt{np_n}}(\varepsilon_{jk}X_{jk})_{j,k=1}^n$.
Let $\lambda_1^{\varepsilon},\ldots,\lambda_n^{\varepsilon}$  denote the (complex)
eigenvalues of the matrix $\bold X^{(\varepsilon)}$  and denote by
$G_n^{\varepsilon}(x,y)$ the empirical spectral distribution
function of the matrix $\bold X^{(\varepsilon)}$,  i. e.
\begin{equation}
\frac1n\sum_{j=1}^nI_{\{\re\{\lambda_j^{\varepsilon}\}\leq x,\;\im\{\lambda_j^{\varepsilon}\}\le y\}}.
\end{equation}
\begin{thm}\label{sparse}
Let $X_{jk}$ be independent identically distributed  sub-Gaussian random variables with
\begin{equation}\notag
\E X_{jk}=0,\qquad \E |X_{jk}|^2=1.
\end{equation}
 Assume that $np_n^{4}\to\infty$ as $n\to\infty$.
Then  $\E G_n^{\varepsilon}(x,y)$ converges weakly to the distribution function $G(x,y)$ as $n\to\infty$.
\end{thm}
\begin{rem}  The assumption $np_n^4\to\infty$ is merely technical and due to our approach to bound the minimal singular values of sparse matrices. For details 
see Subsection \ref{sparse1} in the Appendix.
\end{rem}
\begin{rem}  The crucial problem of the proofs of Theorems \ref{thm0} and \ref{sparse} is
to bound the minimal singular values of shifted matrices $\bold X-z\bold I$ and $\bold X^{\varepsilon}
-z\bold I$. These bounds are based on the results  obtained by Rudelson in \cite{rud:06}.
\end{rem}

The investigation of the convergence the spectral distribution functions
of real or complex (non-symmetric and non-Hermitian) random matrices with
independent entries has a long history. Ginibre in 1965,  \cite{ginibre},
 studied the real, complex and quaternion matrices with i.\ i.\ d.
Gaussian entries. He  derived the joint density for the  distribution of
eigenvalues of matrix. Using the Ginibre results, Edelman in 1997,
\cite{edelman} proved the circular law for the matrices with i.\ i.\ d.
Gaussian entries.   Girko in 1984,
\cite{Girko:1984a}, investigated  the circular law for  general matrices with
independent entries assuming that the distribution of the entries have densities.
  As pointed out by  Bai
\cite{Bai:1997}, Girko's proof had serious gaps. Bai in \cite{Bai:1997}  gave  a  proof of the  circular law
for random matrices with independent entries assuming  that the entries
had bounded 
densities and finite sixth moments.
Unfortunately  this result still does not  cover 
 the case the the Wigner ensemble and  in particular ensembles of  matrices
 with Rademacher entries. These ensembles are of some interest in
various
applications, see e.g.  \cite{Timm:2004}.
  (Wigner,  in his pioneering work in 1955
\cite{Wigner} proved the  {\it semi-circular}  
law for symmetric matrices with i.\ i.\ d.\
Rademacher entries).  A discussion of  Girko's contribution to the proof of
the universality of the  cicular law   may be
found  in Edelman \cite{edelman} as well. Girko published several papers
providing  additional explanations and corrections of  his arguments in his paper in 1984
\cite{Girko:1984a}, see, for example, \cite{Girko:1994}, \cite{Girko:1997},
\cite{Girko:2004}. In \cite{Girko:1994} he states the circular law for
 matrices with independent entries without any assumption on their
densities.  His proof unfortunately does not show  why (assuming his conditions) 
\begin{equation}\notag
\lim_{\varepsilon\to 0}\lim_{n\to\infty}\E\log|\det(\bold X(z)(\bold X(z))^*+\varepsilon^2\bold I|=
\lim_{n\to\infty}\lim_{\varepsilon\to 0}\E\log|\det(\bold X(z)(\bold X(z))^*+\varepsilon^2\bold I|.
\end{equation}
 See for example Khoruzhenko's \cite{Khoruzh:01}, remark on the ``regularization of potential''.
{Girko's  \cite{Girko:1984a} approach  using  families
of spectra of Hermitian matrices for a characterisation of the circular-law
based on the so-called {\it V-transform} was fruitful for all later work. See, for example, Girko's Lemma 1 in \cite{Bai:1997}.}

  We shall outline his approach  using logarithmic potential theory.
Let $\xi$ denote a random variable uniformly distributed over the unit
disc. For any $r>0$, consider the matrix,
$$
\bold X(r)=\bold X-r\xi\bold I,
$$
where $\bold I$ denotes the identity  matrix of order $n$.
Let $\mu_n^{(r)}$ be empirical spectral measure of matrix $\bold X(r)$ defined on the complex plane
as empirical measure of the set of eigenvalues of matrix.
We define a
logarithmic potential of the expected spectral measure $\E\mu_n^{(r)}(ds,dt)$
as
$$
U_n^{(r)}(z)=-\frac1n\E\log |\det(\bold X(r)-z\bold I)|
=-\frac1n\sum\E\log|\lambda_j-z-r\xi|,
$$
where $\lambda_1,\ldots,\lambda_n$ are the eigenvalues of the matrix
$\bold X$. Note that the expected spectral measure $\E\mu_n^{(r)}$ is the
convolution of the measure $\E\mu_n$ and the uniform distribution on the disc
of radius $r$  (see Lemma \ref{ap1} in the Appendix for details).
\begin{lem}\label{smo1} Assume that the sequence $\E\mu_n^{(r)}$ converges weakly  to a measure $\mu$
as $n\to\infty$ and $r\to 0$.
Then
\begin{equation}
\mu=\lim_{n\to\infty}\E\mu_n.
\end{equation}
\end{lem}
\begin{proof}
Let $J$ be a random variable which is uniformly distributed on the set
$\{1,\ldots,n\}$ and independent of the matrix $\bold X$. We may
represent the measure $\E\mu_n^{(r)}$ as distribution of a random variable
$\lambda_J+r\xi$ where $\lambda_J$ and $\xi$ are independent. Computing
the characteristic function of this measure and passing
first to the limit with respect to $n\to\infty$ and then with respect to
$r\to 0$  (see  also Lemma \ref{ap2} in the Appendix), we conclude the result.
\end{proof}
 Now we may fix
$r>0$ and consider the measures $\E\mu_n^{(r)}$. They have bounded
densities. Assume that the measures $\E\mu_n$ have
supports in a fixed compact set  and that $\E\mu_n$ converges weakly to a
measure $\mu$. Applying Theorem 6.9 (Lower Envelope Theorem) from
\cite{saff}, p. 73  (see also Subsection \ref{potential} in the Appendix), we obtain that under these assumptions
\begin{equation}
\liminf_{n\to\infty}U_n^{(r)}(z)=U^{(r)}(z),
\end{equation}
for quasi-everywhere in $\mathbb C$ (for the definition of
``{\it quasi-everywhere}''  see for
example \cite{saff}, p 24  and Subsection \ref{potential} in the Appendix).
Here $U^{(r)}(z)$ denotes the logarithmic potential of measure
$\mu^{(r)}$ which is the convolution of a measure
$\mu$ and of the uniform distribution on the disc of radius $r$.
 Furthermore, note that $U^{(r)}(z)$
we may represented as
\begin{equation}\notag
U^{(r)}(z_0)=\frac2{ r^2}\int_0^rvL(\mu;z_0,v)dv,
\end{equation}
where
\begin{equation}
L(\mu;z_0,v)=\frac1{2\pi}\int_{-\pi}^{\pi}U^{(\mu)}(z_0+v\exp\{i\theta\})d\theta.
\end{equation}
Applying Theorem 1.2 in \cite{saff}, p. 84, (Theorem \ref{thm_saff} in Subsection \ref{potential} in the Appendix)
we get
\begin{equation}\notag
\lim_{r\to0}U_{\mu}^{(r)}(z)=U_{\mu}(z).
\end{equation}
Let $s_1(\bold X)\ge\ldots\ge s_n(\bold X)$ denote the singular values of
 matrix $\bold X$.
Note that for any  $M>2$
\begin{equation}
\Pr\{s_1(\bold X)>M\}\le\Pr\{s_1^2(\bold X)>4\}
\le \sup_x|\E F_n(x)-M_1(x)|\le Cn^{-\frac18},
\end{equation}
where $F_n(x)$ denotes the empirical distribution function of the matrix
 $\frac1n \bold X\bold X^*$.  Here $\bold X^*$ stands for the complex conjugate and transpose of 
the matrix
$\bold X$, and $M_1(x)$
denotes Marchenko--Pastur distribution function with parameter 1 and
density
\begin{equation}\notag
m_1(x)=\frac1{2{\pi}}\sqrt{\frac{4-x}x}I_{\{0< x<4\}}.
\end{equation}
  (See, for example, \cite{Bai:93b}, Theorem 3.2). This implies that
  the  sequence of
measures $\E\mu_n$ is weakly relatively compact. These results imply that we may restrict  the
measures $\E\mu_n$
to some compact set $K$ such that $\sup_n\E\mu_n(K^{(c)})\to 0$. If we take some subsequence of the
sequence of restricted measures $\E\mu_n$ which converges to some measure $\mu$,
then \newline $\liminf_{n\to\infty}U_{\mu_n}^{(r)}(z)=U_{\mu}^{(r)}(z)$, $r>0$ and 
$\lim_{r\to 0}U_{\mu}^{(r)}(z)=U_{\mu}(z)$. If we prove that  $\liminf_{n\to\infty}U_{\mu_n}^{(r)}(z)$ exists
and $U_{\mu}(z)$ is equal to the logarithmic potential
corresponding the uniform distribution on the unit disc then the sequence of
measures $\E\mu_n$ weakly converges to the uniform distribution on the
unit disc. Moreover, it is enough to prove that for some sequence
$r=r(n)\to 0$, $\lim_{n\to\infty}U_{\mu_n}^{(r)}(z)=U_{\mu}(z)$.

Furthermore, let $s_1^{\varepsilon}(z,r)\ge\ldots \ge s_n^{\varepsilon}(z,r)$
denote the singular values of matrix $\bold X^{\varepsilon}(z,r)=\bold
X^{\varepsilon}(r)-z\bold I$. We shall investigate the logarithmic
potential $U_{\mu_n}^{(r)}(z)$.
  Using  elementary properties of  singular values (see for
instance Lemma 3.3 \cite{Goh:69}, p.35), we may represent the function
$U_{\mu_n}^{(r)}(z)$ as follows
\begin{equation}\notag
U_{\mu_n}^{(r)}(z)=-\frac1n\sum_{j=1}^n\E\log{s_j^{\varepsilon}(z,r)}
 =-\frac12\int_0^{\infty}\log
x\E\nu_n^{\varepsilon}(dx,z,r),
\end{equation}
where $\nu_n^{\varepsilon}(\cdot,z,r)$ denotes the spectral measure of the
matrix $\bold H_n^{\varepsilon}(z,r)=(\bold X^{\varepsilon}(r)-z\bold
I)(\bold X^{\varepsilon}(r)-z\bold I)^*$, which is the counting measure of the set of
eigenvalues of the matrix $\bold H_n^{\varepsilon}(z,r)$).

In Section \ref{convergence}) we investigate
 convergence of measure
$\nu_n^{\varepsilon}(\cdot,z)=\nu^{\varepsilon}(\cdot,z,0)$. In Section
\ref{property} we study the properties of the limit measures
$\nu(\cdot,z)$. But the crucial problem for the proof of the circular law is
the so called ``regularization of potential'' problem. See Khoruzhenko
\cite{Khoruzh:01}. We solve this problem using 
 bounds for the minimal singular values
of matrices $\bold X^{\varepsilon}(z):=\bold X^{\varepsilon}-z\bold I$ 
based on techniques developed in Rudelson \cite{rud:06}.
These bounds are given in Section \ref{singular} and in the Appendix, Subsection \ref{sparse1}. 
In Section \ref{proof}
we give the proof of the main Theorem. In the Appendix we combine precise statements of relevant results.
from potential theory and some auxiliary
inequalities for the resolvent matrices.

\section{Convergence of  $\nu_n^{\varepsilon}(\cdot,z)$}\label{convergence}
Denote by $F_n^{\varepsilon}(x,z)$ the distribution function
of the measure $\nu_n^{\varepsilon}(\cdot,z)$,
\begin{equation}\notag
F_n^{\varepsilon}(x,z)=\frac1n\sum_{j=1}^nI_{\{(s_j^{\varepsilon}(z))^2<x\}},
\end{equation}
where $s_1^{\varepsilon}(z)\ge \ldots\ge s_n^{\varepsilon}(z)\ge 0$
denote the singular values of the matrix $\bold X^{\varepsilon}(z)=\bold
X^{\varepsilon}-z\bold I$.  For a
positive random variable $\xi$ and a Rademacher random variable (r. v.)
$\kappa$ consider the transformed r. v.  $\widetilde\xi=\kappa\sqrt{\xi}$. 
 If $\zeta$ has distribution function ${F}_n^{\varepsilon}(x,z)$ the variable $\widetilde{\zeta}$ has distribution 
function $\widetilde{F}_n^{\varepsilon}(x,z)$, given by
$$
\widetilde{F}_n^{\varepsilon}(x,z)=\frac12(1+{\rm{sgn}}
\{x\}F_n^{\varepsilon}(x^2,z))
$$
for all real $x$.
 Note that this induces 
a one-to-one corresponds between the  respective measures
$\nu_n^{\varepsilon}(\cdot,z)$ and ${\widetilde{\nu}}_n^{\varepsilon}(\cdot,z)$. 
The limit distribution function of ${F}_n^{\varepsilon}(x,z)$ as $n\to\infty$ , is denoted
by $F(\cdot,z)$ with corresponding symmetrization $\widetilde{F}(x,z)$
being the limit of $\widetilde{F}_n^{\varepsilon}(x,z)$ as $n\to\infty$. We have
\begin{equation}\notag
\sup_x|F_n^{\varepsilon}(x,z)-F(x,z)|=\sup_x|\widetilde
F_n^{\varepsilon}(x,z)-\widetilde F(x,z)|.
\end{equation}

Denote by $s_n^{\varepsilon}(\alpha,z)$ (resp. $s(\alpha,z)$) and
$S_n^{\varepsilon}(x,z)$ (resp. $S(x,z)$)  the Stieltjes transforms of the
measures $\nu_n^{\varepsilon}(\cdot,z)$ (resp. $\nu(\cdot,z)$) and
$\widetilde {\nu}_n^{\varepsilon}(\cdot,z)$ (resp. $\widetilde
{\nu}(\cdot,z)$) correspondingly. Then we have
\begin{align}\notag
S_n^{\varepsilon}(\alpha, z)=\alpha s_n^{\varepsilon}(\alpha^2, z),\qquad
S(\alpha, z) =\alpha s(\alpha^2, z).
\end{align}
\begin{rem}As is shown in Bai \cite{Bai:1997}, the measure
$\nu(\cdot,z)$ has a density $p(x,z)$
and bounded support. More precisely,
$p(x,z)\le C\max\{1, \frac1{\sqrt x}\}$.
Thus  the measure $\widetilde \nu(\cdot,z)$ has bounded support and
bounded density
$\widetilde p(x,z)=|x|p(x^2,z)$.
\end{rem}
\begin{thm}Let $\E X_{jk}=0$, $\E|X_{jk}|^2=1$, and
\begin{equation}
\varkappa_3=\max_{1\le j,k<\infty}\E|X_{jk}|^3.
\end{equation}
Then
\begin{equation}
\sup_x|F_n^{\varepsilon}(x,z)-F(x,z)|\le C\varkappa_3(np_n)^{-\frac1{10}}.
\end{equation}
\end{thm}
\begin{proof}
To bound the distance between  the distribution functions
 $\widetilde F_n^{\varepsilon}(x,z)$ and $\widetilde F(x,z)$
we investigate the distance between the Stieltjes transforms of these
distribution functions. Introduce the Hermitian $2n\times 2n$ matrix
\begin{equation}\notag
\bold W=\left(\begin{matrix} {\bold O_n\quad(\bold X^{\varepsilon}-z\bold
I)}
\\{(\bold X^{\varepsilon}-z\bold I)^*\quad \bold O_n}\end{matrix}\right),
\end{equation}
where $\bold O_n$ denotes $n\times n$ matrix with all entries equal
to zero.
From \v Sur's complement formula (see for example \cite{HoJohn:91},
Ch. 08, p. 21) it follows that,
for $\alpha=u+iv$, $v>0$,

\begin{equation}\label{shur}
(\bold W-\alpha\bold I_{2n})^{-1}=\left(
\begin{matrix}
{ \alpha\left(\bold X^{\varepsilon}(z)(\bold
X^{\varepsilon}(z))^*-\alpha^2\bold I_n\right)^{-1} \quad \bold
X^{\varepsilon}(z)\left(\bold X^{\varepsilon}(z)(\bold
X^{\varepsilon}(z))^*-\alpha^2\bold I_n\right)^{-1} }
\\
{ \left((\bold X^{\varepsilon}(z))^*\bold
X^{\varepsilon}(z)-\alpha^2\bold I_n\right)^{-1}(\bold
X^{\varepsilon}(z))^* \quad \alpha\left((\bold X^{\varepsilon}(z))^*\bold
X^{\varepsilon}(z)-\alpha^2\bold I_n\right)^{-1} }
\end{matrix}
\right)
\end{equation}
where $\bold X^{\varepsilon}(z)=\bold X^{\varepsilon}-z\bold I$ and
$\bold I_{2n}$ denotes the unit matrix of order $2n$. By definition of
$S_n^{\varepsilon}(\alpha,z)$, we have
\begin{equation}\notag
S_n^{\varepsilon}(\alpha,z)=\frac1{2n}\E\Tr(\bold W-\alpha\bold
I_{2n})^{-1}.
\end{equation}
Set $\bold R(\alpha,z):=(R_{j,k}(\alpha,z))_{j,k=1}^{2n}
=(\bold W-\alpha\bold I_{2n})^{-1}$.
It is easy to check that
\begin{equation}\notag
1+\alpha S_n^{\varepsilon}(\alpha,z)=\frac{1}{2n}\E\Tr\bold W\bold
R(\alpha,z).
\end{equation}
We may rewrite this equality  as
\begin{align}\label{n1}
1+\alpha S_n^{\varepsilon}(\alpha,z)&=\frac1{2n\sqrt
{np_n}}\sum_{j,k=1}^n \E
({\varepsilon_{jk}}X_{jk}R_{k+n,j}(\alpha,z)+{\varepsilon_{jk}}\overline
X_{jk}
R_{j+n,k}(\alpha,z))\notag\\
&-\frac{\overline z}{2n}\sum_{j=1}^n\E R_{j,j+n}(\alpha,z)
-\frac{z}{2n}\sum_{j=1}^n\E R_{j+n,j}(\alpha,z).
\end{align}
We introduce the notations
\begin{align}\notag
\bold A&=(\bold X^{\varepsilon}(z)(\bold
X^{\varepsilon}(z))^*-\alpha^2\bold I)^{-1},
\quad \bold B=\bold X^{\varepsilon}(z)\bold A,\notag\\
\bold C&=((\bold X^{\varepsilon}(z))^*\bold
X^{\varepsilon}(z)-\alpha^2\bold I)^{-1}, \quad \bold D=\bold C(\bold
X^{\varepsilon}(z))^*.\notag
\end{align}
With these notations we rewrite equality (\ref{shur}) as follows
\begin{equation}\label{shur1}
\bold R(\alpha,z)=(\bold W-\alpha\bold I_{2n})^{-1}=\left(
\begin{matrix}
{
\alpha\bold A
\quad
\bold B}
\\
{
\bold D\quad
\alpha\bold C}
\end{matrix}
\right)
\end{equation}
Equalities (\ref{shur1}) and (\ref{n1}) together imply
\begin{align}\label{n2}
1+\alpha S_n^{\varepsilon}(\alpha,z)&=\frac1{2n\sqrt {np_n}}\sum_{j,k=1}^n
\E ({\varepsilon_{jk}}X_{jk}R_{k+n,j}(\alpha,z)+\varepsilon_{jk}\overline X_{jk}R_{j,k+n}(\alpha,z))\notag\\
&-\frac{z}{2n}\E \Tr \bold D-\frac{\overline z}{2n}\E \Tr \bold B.
\end{align}

In the what follows we shall use a simple resolvent equality.
For two matrices
$\bold U$ and $\bold V$
let $\bold R_{U}=(\bold U-\alpha\bold I)^{-1}$, $\bold R_{U+V}=
(\bold U+\bold V-\alpha\bold I)^{-1}$,
then
\begin{equation}\notag
\bold R_{U+V}=\bold R_{U}-\bold R_{U}\bold V\bold R_{U+V}.
\end{equation}
Let $\{\bold e_1,\ldots \bold e_{2n}\}$ denote the canonical orthonormal basis
in $\mathbb R^{2n}$.
Let $\bold W^{(jk)}$ denote the matrix is obtained from
$\bold W$ by replacing the both entries $X_{j,k}$ and
$\overline X_{j,k}$ by 0. In our notation we may write
\begin{equation}\label{repr1}
\bold W=\bold W^{(jk)}+\frac1{\sqrt {np_n}}{\varepsilon_{jk}}X_{jk}\bold
e_{j}\bold e_{k+n}^T+ \frac1{\sqrt {np_n}}{\varepsilon_{jk}}\overline
X_{jk}\bold e_{k+n}\bold e_{j}^T.
\end{equation}
Using this representation and the resolvent equality, we get
\begin{equation}\label{rezolvent}
\bold R=\bold R^{(j,k)}-\frac1{\sqrt {np_n}}{\varepsilon_{jk}}X_{jk}\bold
R^{(j,k)}\bold e_{j} \bold e_{k+n}^T\bold R -\frac1{\sqrt
{np_n}}{\varepsilon_{jk}}\overline X_{jk}\bold R^{(j,k)}\bold
e_{k+n}\bold e_{j}^T \bold R.
\end{equation}
Here and in the what follows we omit the arguments $\alpha$ and $z$ in
the notation of resolvent matrices. For any vector $\bold a$, let $\bold a^T$
denote the transposed vector $\bold a$. Applying the resolvent equality again,
we obtain
\begin{align}\notag
\bold R=\bold R^{(j,k)}-\frac1{\sqrt {np_n}}{\varepsilon_{jk}}X_{jk}\bold
R^{(j,k)}\bold e_{j} \bold e_{k+n}^T\bold R^{(j,k)} -\frac1{\sqrt
{np_n}}{\varepsilon_{jk}}\overline X_{jk}\bold R^{(j,k)}\bold
e_{k+n}\bold e_{j}^T \bold R^{(j,k)} +\bold T^{(jk)},
\end{align}
where
\begin{align}\notag
\bold T^{(jk)}&=\frac1{{np_n}}{\varepsilon_{jk}}X_{jk}^2\bold
R^{(j,k)}\bold e_{j}\bold e_{k+n}^T \bold R^{(j,k)}\bold e_{j}
\bold e_{k+n}^T\bold R\notag\\
&+\frac1{{np_n}}{\varepsilon_{jk}}|X_{jk}|^2\bold R^{(j,k)}\bold
e_{j}\bold e_{k+n}^T
\bold R^{(j,k)}\bold e_{k+n}\bold e_{j}^T\bold R\notag\\
&+\frac1{{np_n}}{\varepsilon_{jk}}(\overline X_{jk})^2\bold R^{(j,k)}\bold e_{k+n}\bold e_{j}^T\bold R^{(j,k)}\bold e_{j}\bold e_{k+n}^T\bold R\notag\\
&+\frac1{{np_n}}{\varepsilon_{jk}}|X_{jk}|^2\bold R^{(j,k)}\bold
e_{k+n}\bold e_{j}^T\bold R^{(j,k)}\bold e_{k+n}\bold e_{j}^T\bold
R\notag
\end{align}
This implies
\begin{align}\notag
\bold R_{j,k+n}&=\bold R^{(j,k)}_{j,k+n}-\frac1{\sqrt
{np_n}}{\varepsilon_{jk}}X_{jk} \bold R^{(j,k)}_{j,j}\bold
R^{(j,k)}_{k+n,k+n} -\frac1{\sqrt {np_n}}{\varepsilon_{jk}}\overline
X_{jk}(\bold R^{(j,k)}_{j,k+n})^2+
\bold T^{(j,k)}_{j,k+n}\notag\\
\bold R_{k+n,j}&=\bold R^{(j,k)}_{k+n,j}-\frac1{\sqrt
{np_n}}{\varepsilon_{jk}}X_{jk} \bold R^{(j,k)}_{k+n,j}\bold
R^{(j,k)}_{j,k+n} -\frac1{\sqrt {np_n}}{\varepsilon_{jk}}\overline
X_{jk}\bold R^{(j,k)}_{k+n,k+n} \bold R^{(j,k)}_{j,j}+\bold
T^{(j,k)}_{k+n,j}.
\end{align}
Applying these notations to the equality (\ref{n2}) and taking into
account that $X_{jk}$ and
$\bold R^{(jk)}$ are independent, we get
\begin{align}\label{new12}
1+\alpha S_n^{\varepsilon}(\alpha,z)+\frac{z}{2n}\Tr \bold D+\frac{\overline z}{2n}
\Tr\bold B&=-\frac1{n^2p_n}\sum_{j,k=1}^n \E
{\varepsilon_{jk}}R^{(j,k)}_{j,j} R^{(j,k)}_{k+n,k+n}\notag\\&-
\frac1{2n^2p_n}\sum_{j,k=1}^n \E {\varepsilon_{jk}}|X_{jk}|^2\E
(R^{(j,k)}_{j,k+n})^2\notag
\\&-\frac1{2n\sqrt {np_n}}\sum_{j,k=1}^n
\E ({\varepsilon_{jk}}X_{jk}T^{(j,k)}_{k+n,j}+{\varepsilon_{jk}}\overline
X_{jk}T^{(j,k)}_{j,k+n}) .
\end{align}
By definition of $\bold T^{(j,k)}$ and standard resolvent properties,
we obtain the following bounds,
for any $p,q=1,\ldots,2n$, $j,k=1,\ldots n$,  and any $z=u+iv,\ v>0$,
\begin{align}\label{bound1}
|R_{p,p}-R^{(j,k)}_{pp}|&\le \frac {C{\varepsilon_{jk}}|X_{jk}|}
{\sqrt {np_n}}(|R^{jk}_{pj}||R_{k+n,p}|+|R^{jk}_{p,k+n}||R_{jp}|)\notag\\
\frac1{n^2}\sum_{j,k=1}^n\E |R^{(j,k)}_{j,k+n}|^2&\le
\frac{C}{nv^4}\\
\frac1{n\sqrt
{np_n}}\sum_{j,k=1}^n\E{\varepsilon_{jk}}|X_{jk}||T^{(j,k)}_{j,k+n}|&\le
\frac{C\varkappa_3}{{np_n}v^4}
\end{align}
For the proof of these inequalities see in the Appendix, Lemma \ref{ap}.
Using the last inequalities we obtain, that for $v>0$
\begin{align}\label{41}
\left|\frac1n\sum_{j=1}^n\E
R_{jj}\frac1n\sum_{k=1}^nR_{k+n,k+n}\right.&\left.-\frac1{n^2}
\sum_{j=1}^n\sum_{k=1}^n\E R^{(jk)}_{jj}
R^{(jk)}_{k+n,k+n}\right|\notag\\&\le \frac{C}{n^2\sqrt
{np_n}v}\sum_{j=1}^n \sum_{k=1}^n\E\varepsilon_{jk}|X_{jk}|
(|R^{(jk)}_{jj}||R_{k+n,j}| +
|R^{(jk)}_{j,k+n}||R_{jj}|)\notag\\
&\le \frac{C\varkappa_3}{nv^4}.
\end{align}

Since $\frac1n\sum_{j=1}^nR_{jj}=\frac1n\sum_{k=1}^nR_{k+n,k+n}
=\frac1{2n}\Tr \bold R(\alpha,z)$,
we obtain
\begin{equation}\label{prev}
|\frac1{n^2}\sum_{j=1}^n\sum_{k=1}^n\E R^{(jk)}_{jj}
R^{(jk)}_{k+n,k+n}-\E(\frac1{2n}\Tr \bold R(\alpha,z))^2|\le \frac
{C\varkappa_3}{nv^4}
\end{equation}
Note that for any Hermitian random matrix $\bold W$ with independent
entries on and above the
diagonal we have
\begin{equation}\label{variation}
\E|\frac 1n\Tr \bold R(\alpha,z)-\E\frac 1n\Tr \bold R(\alpha,z)|^2\le
\frac{C}{nv^2}.
\end{equation}
The proof of this inequality is easy and due to a martingale type 
expansion already used by Girko. Inequalities (\ref{prev}) and (\ref{variation}) together imply
that for $v>cn^{-\frac14}$
\begin{equation}\label{ne6}
|\frac1{n^2}\sum_{j=1}^n\sum_{k=1}^n\E R^{(jk)}_{jj}
R^{(jk)}_{k+n,k+n}-(S_n^{\varepsilon}(\alpha,z))^2|\le \frac{C}{ \sqrt nv^2}
\end{equation}

We may now rewrite equality (\ref{n2}) as follows
\begin{align}\label{n3}
1+\alpha S_n^{\varepsilon}(\alpha,z)+(S_n^{\varepsilon}(\alpha,z))^2= -\frac{z}{2n}\E \Tr \bold
D-\frac{\overline z}{2n}\E \Tr \bold B+ \theta\frac1{ \sqrt{np_n}v^2},
\end{align}
were $\theta$ is a function such that $|\theta|\le 1$ and
$v>c(np_n)^{-\frac14}$.

 We now investigate the functions
$T(\alpha,z)=\frac1n\E \Tr \bold D$ and $V(\alpha,z)=\frac1n\E\bold B$.
Since the arguments for  both functions are similar we provide  it for the
first one only. By definition of the matrix $\bold B$, we have
\begin{equation}\notag
\Tr \bold B=\frac1{\sqrt
{np_n}}\sum_{j,k=1}^n\varepsilon_{jk}X_{j,k}((\bold
X^{\varepsilon}(z)(\bold X^{\varepsilon}(z))^* -\alpha^2)^{-1})_{kj}
\end{equation}
According to equality (\ref{shur1}), we have
\begin{equation}\notag
\Tr \bold B=\frac1{\alpha\sqrt
{np_n}}\sum_{j,k=1}^n\varepsilon_{jk}X_{j,k}R_{kj}-z \Tr \bold A
\end{equation}

Using the resolvent equality (\ref{rezolvent}) and Lemma \ref{ap}, we get,
for $v>c(np_n)^{-\frac14}$
\begin{align}\label{n7}
T(\alpha,z)=-\frac1{\alpha n^2}\sum_{j,k=1}^n\E R^{(jk)}_{k,k+n}
R^{(jk)}_{jj}-\frac{z}{\alpha}S_n^{\varepsilon}(\alpha,z) +\theta\frac{C\varkappa_3}{
{np_n}v^2}.
\end{align}

Similar to (\ref{ne6}) we obtain
\begin{equation}\label{n8}
|\frac1{\alpha n^2}\sum_{j,k=1}^n\E R^{(jk)}_{jj}R^{(jk)}_{k,k+n}-
V(\alpha,z)S_n^{\varepsilon}(\alpha,z)|\le \frac C{n\sqrt nv^4}
\end{equation}
Inequalities (\ref{n7}) and (\ref{n8}) together imply, for
$v>c(np_n)^{-\frac14}$,
\begin{equation}\label{ne8}
V(\alpha,z)=-\frac {\overline zS_n^{\varepsilon}(\alpha,z)}{\alpha+S_n^{\varepsilon}(\alpha,z)}
+\theta\frac{C\varkappa_3}{  {np_n} v^2|\alpha+S_n^{\varepsilon}(\alpha,z)|}.
\end{equation}
Analogously we get
\begin{equation}\label{ne9}
T(\alpha,z)=-\frac {zS_n^{\varepsilon}(\alpha,z)}{\alpha+S_n^{\varepsilon}(\alpha,z)}
+\theta\frac{C}{ {np_n} v^2|\alpha+S_n^{\varepsilon}(\alpha,z)|}.
\end{equation}
Insecting (\ref{ne8}) and (\ref{ne9}) in (\ref{new12}), we get
\begin{equation}\label{ne10}
(S_n^{\varepsilon}(\alpha,z))^2+\alpha S_n^{\varepsilon}(\alpha,z)+1-\frac{|z|^2S_n^{\varepsilon}(\alpha,z)}
{\alpha+S_n^{\varepsilon}(\alpha,z)}=\delta_n(z),
\end{equation}
where
\begin{equation}\notag
|\delta_n(\alpha,z)|\le \frac{C\varkappa_3}{ {np_n}
v^2|S_n^{\varepsilon}(\alpha,z)+\alpha|}.
\end{equation}
or equivalently
\begin{equation}\label{99}
S_n^{\varepsilon}(\alpha,z)\left(\alpha+S_n^{\varepsilon}(\alpha,z)\right)^2
+\left(\alpha+S_n(^{\varepsilon}\alpha,z)\right)-|z|^2S_n^{\varepsilon}(\alpha,z)=
\widetilde\delta_n(\alpha,z),
\end{equation}
were $\widetilde\delta_n(\alpha,z)=\theta\frac{C\varkappa_3}{ {np_n}v^2}$.
The last equation we may rewrite as
\begin{equation}\label{95}
S_n^{\varepsilon}(\alpha,z)=-\frac{\alpha+S_n^{\varepsilon}(\alpha,z)}
{(\alpha+S_n^{\varepsilon}(\alpha,z))^2-|z|^2}+\widehat\delta_n(\alpha,z),
\end{equation}
were
\begin{equation}\label{delta}
\widehat\delta_n(\alpha,z)=\frac{\widetilde\delta_n(\alpha,z)}
{(\alpha+S_n^{\varepsilon}(\alpha,z))^2-|z|^2}.
\end{equation}
Note that
\begin{equation}\notag
\left|\frac{1}{(\alpha+S_n^{\varepsilon}(\alpha,z))^2-|z|^2}\right| \le
\frac1{v|\alpha+S_n^{\varepsilon}(\alpha,z)|}
\end{equation}.
This implies that
\begin{equation}\notag
|\widehat\delta_n(\alpha,z)| \le \frac C{ np_nv^2|\alpha+S_n^{\varepsilon}(\alpha,z)|}.
\end{equation}


Furthermore, we prove the following simple Lemma.
\begin{lem}\label{lem03}Let $\alpha=u+iv$, $v>0$. Let $S(\alpha, z)$ satisfy the equation
\begin{equation}\label{eq:1}
S(\alpha,z)=-\frac{\alpha+S(\alpha,z)}{(\alpha+S(\alpha,z))^2-|z|^2}.
\end{equation}
and $\im\{S(\alpha,z)\}>0$.
Then the following inequality
\begin{equation}\notag
1-|S(\alpha,z)|^2-\frac{|z|^2|S(\alpha,z)|^2}{|\alpha+S(\alpha,z)|^2}
\ge\frac v{v+1}.
\end{equation}
holds.
\end{lem}
\begin{proof}
The Stieltjes transform $S(\alpha, z)$ satisfies the following equation,
 for $\alpha=u+iv$ with $v>0$,
\begin{equation}\label{eq:100}
S(\alpha,z)=-\frac{\alpha+S(\alpha,z)}{(\alpha+S(\alpha,z))^2-|z|^2}.
\end{equation}
Comparing the imaginary parts of both sides of this equation, we get
\begin{equation}\label{eq:2}
\im\{\alpha+S(\alpha,z)\}=\im\{\alpha+S(\alpha,z)\}\,\frac{|\alpha
+S(\alpha,z)|^2+|z|^2}{|(\alpha+S(\alpha,z))^2-|z|^2|^2} +v.
\end{equation}
Equations (\ref{eq:1}) and (\ref{eq:2}) together imply
\begin{equation}\label{eq:3}
\im\{\alpha+S(\alpha,z)\}\left(1-\frac{|\alpha
+S(\alpha,z)|^2+|z|^2}{|(\alpha+S(\alpha,z))^2-|z|^2|^2}\right) =v.
\end{equation}
Since $v>0$ and $\im\{\alpha+S(\alpha,z)\}>0$, it follows that
$$
1-\frac{|\alpha
+S(\alpha,z)|^2+|z|^2}{|(\alpha+S(\alpha,z))^2-|z|^2|^2}>0.
$$
In particular, we have
$$
|S(\alpha,z)|\le 1.
$$
Inequality (\ref{eq:3}) and  the last remark together imply
\begin{equation}\notag
1-\frac{|\alpha +S(\alpha,z)|^2+|z|^2}{|(\alpha+S(\alpha,z))^2-|z|^2|^2}
=\frac v{\im\{\alpha+S(\alpha,z)\}}\ge \frac v{v+1}.
\end{equation}
The proof  is completed.
\end{proof}
To compare the function $S(\alpha,z)$ and $S_n(\alpha,z)$ we prove
\begin{lem}\label{lem030}
Let
$$
|\widehat\delta_n(\alpha,z)|\le \frac v2.
$$
Then the following inequality holds
\begin{equation}\notag
1-\frac{|\alpha
+S_n^{\varepsilon}(\alpha,z)|^2+|z|^2}{|(\alpha+S_n^{\varepsilon}(\alpha,z))^2-|z|^2|^2}\ge \frac
v{4}.
\end{equation}
\end{lem}
\begin{proof}
By assumption, we have
$$
\im\{\widehat\delta_n(\alpha, z)+\alpha\}>\frac v2.
$$
Repeating the arguments of Lemma \ref{lem03} completes the proof.
\end{proof}
The next Lemma give as a bound for the distance between the Stieltjes
transforms $S(\alpha,z)$ and $S_n^{\varepsilon}(\alpha, z)$.
\begin{lem}\label{lem04}
Let
$$
|\widehat\delta_n(\alpha,z)|\le \frac v8.
$$
Then
$$
|S_n^{\varepsilon}(\alpha,z)-S(\alpha,z)|\le \frac{4|\widehat\delta_n(\alpha,z)|}v.
$$
\end{lem}
\begin{proof} Note that $S(\alpha,z)$ and $S_n^{\varepsilon}(\alpha,z)$ satisfy the
equations
\begin{equation}\label{eq:10}
S(\alpha,z)=-\frac{\alpha+S(\alpha,z)}{(\alpha+S(\alpha,z)^2-|z|^2}
\end{equation}
and
\begin{equation}\label{eq:110}
S_n^{\varepsilon}(\alpha,z)=-\frac{\alpha+S_n^{\varepsilon}(\alpha,z)}{(\alpha+S_n^{\varepsilon}(\alpha,z)^2
-|z|^2}+\widehat\delta_n(\alpha,z)
\end{equation}
respectively.
These equations together imply
\begin{align}
S(\alpha,z)-S_n^{\varepsilon}(\alpha,z)=\frac{(\alpha+S_n^{\varepsilon}(\alpha,z))(\alpha
+S(\alpha,z))+|z|^2}
{((\alpha+S(\alpha,z)^2-|z|^2)((\alpha+S_n^{\varepsilon}(\alpha,z)^2-|z|^2)}\notag\\
\times(S(\alpha,z)-S_n^{\varepsilon}(\alpha,z))+\widehat\delta_n(\alpha,z).
\end{align}
Applying inequality $|ab|\le \frac12{(a^2+b^2)}$, we get
\begin{align}\notag
\left|1-\frac{(\alpha+S_n^{\varepsilon}(\alpha,z))(\alpha +S(\alpha,z))+|z|^2}
{((\alpha+S(\alpha,z)^2-|z|^2)((\alpha+S_n^{\varepsilon}(\alpha,z)^2-|z|^2)}\right|&
\notag
\\\ge \frac12\left(1-\frac{|\alpha
+S_n(\alpha,z)|^2+|z|^2}{|(\alpha+S_n^{\varepsilon}(\alpha,z))^2-|z|^2|^2}\right)&\notag\\+
\frac12\left(1-\frac{|\alpha
+S(\alpha,z)|^2+|z|^2}{|(\alpha+S(\alpha,z))^2-|z|^2|^2}\right).&\notag
\end{align}
The last inequality and Lemmas \ref{lem03} and \ref{lem030}
together imply
\begin{equation}\notag
\left|1-\frac{(\alpha+S_n^{\varepsilon}(\alpha,z))(\alpha +S(\alpha,z))+|z|^2}
{((\alpha+S(\alpha,z)^2-|z|^2)((\alpha+S_n^{\varepsilon}(\alpha,z)^2-|z|^2)}\right| \ge
\frac v4.
\end{equation}
This completes the proof of the Lemma.
\end{proof}

To bound the distance between the distribution function $F_n(x,z)$ and the 
distribution function $F(x,z)$
corresponding the Stieltjes transform $S(\alpha,z)$ we use Corollary 2.3
from \cite{GT03}.
In the next lemma we give an integral bound for the distance between the
Stieltjes transforms $S(\alpha,z)$ and $S_n^{\varepsilon}(\alpha,z)$.
\begin{lem}
For $v\ge v_0(n)=c(np_n)^{-1/4}$ the inequality
\begin{equation}\notag
\int_{-\infty}^{\infty}|S(\alpha,z)-S_n^{\varepsilon}(\alpha,z)|du\le
\frac{C(1+|z|^2)\varkappa_3}{np_n v^6}.
\end{equation}
holds.
\end{lem}
\begin{proof}
It is enough to prove that

\begin{equation}\notag
\int_{-\infty}^{\infty}|\widehat\delta_n(\alpha, z)|du\le C\gamma_n,
\end{equation}
where $\gamma_n=\frac{C}{np_n v^5}$. By definition of
$\widehat\delta(\alpha,z)$, we have
\begin{align}\label{np3}
\int_{-\infty}^{\infty}|\widehat\delta_n(\alpha, z)|du\le
\frac{c\varkappa_3}{np_nv^2}\int_{-\infty}^{\infty}
\frac{du}{|(\alpha+S_n^{\varepsilon}(\alpha,z))^2-|z|^2|}.
\end{align}
Furthermore,  the representation (\ref{95}) implies that
\begin{align}
\frac{1}{|(\alpha+S_n^{\varepsilon}(\alpha,z))^2-|z|^2|}\le\frac{|S_n^{\varepsilon}(\alpha,z)|}{|\alpha+S_n^{\varepsilon}(\alpha,z)|}
+\frac{|\widehat\delta_n(\alpha,z)|}{|\alpha+S_n^{\varepsilon}(\alpha,z)|}.
\end{align}
Note that, according  to the relation (\ref{99}),
\begin{equation}\label{012}
\frac1{|\alpha+S_n^{\varepsilon}(\alpha,z)|}\le \frac{|z|^2|S_n^{\varepsilon}(\alpha,z)|}
{|\alpha+S_n^{\varepsilon}(\alpha,z)|^2}+|S_n^{\varepsilon}(\alpha,z)|+
\frac{|\widetilde\delta_n(\alpha,z)|}{|\alpha+S_n^{\varepsilon}(\alpha,z)|^2}\le
|S_n^{\varepsilon}(\alpha,z)|(1+\frac{|z|^2}{v^2})+\frac{|\delta_n(\alpha,z)|}
{|\alpha+S_n^{\varepsilon}(\alpha,z)|}.
\end{equation}
This inequality implies
\begin{equation}\label{007}
\int_{-\infty}^{\infty}\frac{|S_n^{\varepsilon}(\alpha,z)|}{|\alpha+S_n^{\varepsilon}(\alpha,z)|}du\le
\frac {C(1+|z|^2)}{v^2}\int_{-\infty}^{\infty}|S_n^{\varepsilon}(\alpha,z)|^2du+
\int_{-\infty}^{\infty}|\delta_n(\alpha,z)|
\frac{|S_n^{\varepsilon}(\alpha,z)|}{|\alpha+S_n^{\varepsilon}(\alpha,z)|}du
\end{equation}
 It follows  from the relation (\ref{ne10}), for $v>c(np_n)^{-\frac14}$, that
\begin{equation}
|\delta_n(\alpha,z)|\le \frac{C\varkappa_3}{np_nv^3}.
\end{equation}
The last two inequalities together imply that for sufficiently large $n$
and $v>c(np_n)^{-\frac14}$,
\begin{equation}
\int_{-\infty}^{\infty}\frac{|S_n^{\varepsilon}(\alpha,z)|}{|\alpha+S_n^{\varepsilon}(\alpha,z)|}du\le
 \frac
{C(1+|z|^2)}{v^2}\int_{-\infty}^{\infty}|S_n^{\varepsilon}(\alpha,z)|^2du
\le\frac{C(1+|z|^2)}{v^3} .
\end{equation}
 The inequalities
(\ref{012}), (\ref{np3}), and the definition of
$\widehat\delta_n(\alpha,z)$ together imply
\begin{equation}
\int_{-\infty}^{\infty}|\widehat\delta_n(\alpha,
z)|du\le\frac{C(1+|z|^2)}{ np_nv^5}
+\frac{C\varkappa_3}{np_nv^3}\int_{-\infty}^{\infty}|\widehat\delta_n(\alpha,
z)|du.
\end{equation}
If we choose $v$ such that $\frac C{ np_nv^3}<\frac12$ we obtain
\begin{equation}
\int_{-\infty}^{\infty}|\widehat\delta_n(\alpha, z)|du\le
\frac{C(1+|z|^2)}{np_nv^5} .
\end{equation}

\end{proof}

In Section \ref{property} is shown that the measure $\nu(\cdot,z)$ has
bounded support and bounded density for any $z$. To bound the distance
between  the distribution functions $\E F_n(x,z)$ and $F(x,z)$ we may apply
Corollary 3.2 from \cite{GT03}  (see also Lemma \ref{ap3} in the Appendix). We take $V=1$ and
$v_0=C(np_n)^{-\frac14}$. Then Lemmas \ref{lem03} and \ref{lem030}
together imply

\begin{equation}\label{supremum}
\sup_x|\E F_n^{\varepsilon}(x,z)-F(x,z)|\le C(np_n)^{-\frac1{4}}.
\end{equation}
\end{proof}
\section{Properties of the measure $\nu(\cdot,z)$}\label{property}
In this Section we investigate the properties of the measure $\nu(\cdot,z)$. At
first note that there \mc{exists} a solution $S(\alpha,z)$ of the equation
\begin{equation}\label{77}
S(\alpha,z)=-\frac{S(\alpha,z)+z}{(S(\alpha,z)+z)^2-|z|^2}
\end{equation}
such that
$$
\im\{S(\alpha,z)\mc{\}}\ge 0\quad\text{for}\quad v>0
$$
and $S(\alpha,z)$ is an analytic function in the upper half-plane
$\alpha=u+iv$,  $v>0$.  This follows from the relative compactness of the sequence of
analytic functions $S_n(\alpha,z)$, $n\in\mathbb N$. From  (\ref{eq:10}) it follows
immediately that
\begin{equation}
|S(\alpha,z)|\le 1.
\end{equation}
Set $y=S(x,z)+x$ and consider the equation (\ref{eq:10}) on the real line
\begin{equation}
y=-\frac{y}{y^2-|z|^2}+x,
\end{equation}
or
\begin{equation}\label{eq:11}
y^3-xy^2+(1-|z|^2)y+x|z|^2=0.
\end{equation}

Set
\begin{equation}
x_1^2=\frac{5+2|z|^2}{2}+\frac{(1+8|z|^2)^{\frac32}-1}{8|z|^2 },\quad
x_2^2=\frac{5+2|z|^2}{2}-\frac{(1+8|z|^2)^{\frac32}+1}{8|z|^2 }.
\end{equation}
It is straightforward to check that for $|z|\le 1$
$\sqrt{3(1-|z|^2)}\le|x_1|$
and
$x_2^2<0$ for $|z|<1$ and
$x_2^2=0$ for $|z|=1$, and $x_2^2>0$ for $|z|>1$.
\begin{lem}\label{lem01}
In the case $|z|\le 1$ equation (\ref{eq:11}) has one real root for
$|x|\le|x_1|$
and three real roots for
$|x|>|x_1|$.
In the case $|z|>1$ equation (\ref{eq:11}) has one real root for
$|x_2|\le x\le |x_1|$
and has tree real roots for
$|x|\le |x_2|$
or for
$|x|\ge |x_1|$.
\end{lem}

\begin{proof}
Set
$$
L(y):=y^3-xy^2+(1-|z|^2)y+x|z|^2.
$$
We consider the roots equation
\begin{equation}
L'(y)=3y^2-2xy+(1-|z|^2)=0.
\end{equation}
The roots of this equation are
$$
y_{1,2}=\frac{x\pm\sqrt{x^2-3(1-|z|^2)}}3.
$$
This implies that, for $|z|\le 1$ and for
$$
|x|\le \sqrt{3(1-|z|^2)},
$$
the equation (\ref{eq:11}) has one real root.
Furthermore, direct calculations shown that
$$
L(y_1)L(y_2)=\frac1{27}\left(-4|z|^2x^4+(8|z|^4+20|z|^2-1)x^2
+4(1-|z|^2)^3\right)
$$

Solving the  equation $L(y_1)L(y_2)=0$ with respect to $x$, we get
for $|z|\le 1$ and $\sqrt{3(1-|z|^2)}\le |x|\le|x_1|$
$$
L(y_1)L(y_2)\ge 0,
$$
and
for $|z|\le1$ and $|x|>\sqrt{\frac{20+8|z|^2}{8}
+\frac{(1+8|z|^2)^{\frac32}-1}
{8|z|^2 }}$
$$
L(y_1)L(y_2)< 0,
$$
These relations imply that
for $|z|\le 1$ the function
$L(y)$ has three real roots for $|x|\ge|x_1|$
and one real root for $|x|<|x_1|$.

Consider the case $|z|>1$ now. In this case $y_{1,2}$ are real for
 all $x$ and $x_2^2>0$.
Note that
$$
L(y_1)L(y_2)\le0
$$
for $|x|\le|x_2|$ and for $|x|\ge |x_1|$ and
$$
L(y_1)L(y_2)>0
$$
for $|x_2|<x<|x_1|$. These implies that for $|z|>1$ and for
$|x_2|<x<|x_1|$ the function $L(y)$ has one real root and for
$|x|\le|x_2|$ or for $|x|\ge |x_1|$ the function $L(y)$ has three real
roots. The Lemma is proved.
\end{proof}
\begin{rem}\label{rem:01}
From Lemma \ref{lem01} it follows that the 
measure $\nu(x,z)$ has
a density $p(x,z)$ and 
 \begin{itemize}
\item{$p(x,z)\le 1$, for all $x$ and $z$}
\item{for $|z|\le 1$, if $|x|\ge x_1$ then $p(x,z)=0$;}
\item{for $|z|\ge 1$, if $|x|\ge x_1$ or $|x|\le x_2$ then $p(x,z)=0$;}
\item{$p(x,z)>0$ otherwise.}
\end{itemize}
\end{rem}
The next lemma is  an analogue of Lemma 4.4 in Bai \cite{Bai:1997}.
\begin{lem}
The following equality
\begin{equation}
\frac{ \partial}{ \partial s}\left(\int_0^{\infty}\log x\nu(dx,z)\right)
=\frac12\Re\{g(x,z)\}
\end{equation}
holds.
\end{lem}
\begin{proof}Following Bai \cite{Bai:1997} Lemma 4.4, we consider
\begin{equation}
I(C):=\int_0^C\frac{\partial y(x)}{\partial s}dx.
\end{equation}
We have
\begin{equation}\label{eq:17}
y^3+2xy^2+x^2y-|z|^2y+y+x=0.
\end{equation}
Taking the derivatives with respect to $x$ and $s$ correspondingly,
we get
\begin{equation}
\frac{\partial y}{\partial x}\left(3y^2+4xy+(1-|z|^2+x^2)\right)
=-1-2y(x+y)
\end{equation}
and
\begin{equation}
\frac{\partial y}{\partial s}\left(3y^2+4xy+(1-|z|^2+x^2)\right)=2sy.
\end{equation}
These equalities together imply
\begin{equation}
\frac{\partial y}{\partial s}=-\frac{2sy}{1+2y(x+y)}\frac{\partial y}
{\partial x}.
\end{equation}
From equation (\ref{eq:17}) it follows that
\begin{equation}\label{eq:18}
1+2y(y+x)=\pm\sqrt{1+4|z|^2y^2}.
\end{equation}
 Using the results of Remark \ref{rem:01}, it is straightforward to check that
for $|z|\le 1$
\begin{equation}
1+2y(y+x)=\sqrt{1+4|z|^2y^2}
\end{equation}
and for $|z|>1$ there exists a number $x_0$ such that
$\sqrt{1+4|z|^2y^2}=0$.  Furthermore, we have for $-x_0\le x\le 0$
\begin{equation}
1+2y(y+x)=\sqrt{1+4|z|^2y^2}
\end{equation}
and
for $x<-x_0$ we obtain
\begin{equation}
1+2y(y+x)=-\sqrt{1+4|z|^2y^2}.
\end{equation}
Using these equalities, we get
\begin{equation}
\int_{-C}^{0}\frac{\partial y}{\partial s}dx=-\int_{-C}^{0}
\frac{2sy}{1+2y(x+y)}\frac{\partial y}{\partial x}dx.
\end{equation}
For $|z|\le 1$, we have
\begin{equation}
\int_{-C}^0\frac{\partial y}{\partial s}dx=-\int_{-C}^{0}
\frac{2sy}{\sqrt{1+4|z|^2y^2}}
\frac{\partial y}{\partial x}dx=\frac{s}{4|z|^2}
\left(\sqrt{1+4|z|^2y^2(-C)}+\sqrt{1+4|z|^2(|z|^2-1)}\right).
\end{equation}
In the limit $C\to\infty$, we get, for $|z|\le 1$,
\begin{equation}\label{rel01}
\int_{-\infty}^0\frac{\partial y}{\partial s}dx=\frac s2.
\end{equation}
For $|z|>1$, we have
\begin{equation}\label{rel02}
\int_{-\infty}^0\frac{\partial y}{\partial s}dx=\int_{-x_0}^0\frac{2sy}
{\sqrt{1+4|z|^2y^2}}
\frac{\partial y}{\partial x}dx-\int_{-\infty}^{-x_0}\frac{2sy}
{\sqrt{1+4|z|^2y^2}}
\frac{\partial y}{\partial x}dx=\frac s{2|z|^2}.
\end{equation}
Similar to Bai \cite{Bai:1997} (equality (4.39)) we have
\begin{align}\label{rel03}
\int_{-C}^0y(x)dx&=\int_{-C}^0y(x)dx=\int_{0}^C\int_0^{\infty}
\frac1{u+x}\nu(du,z)dx\notag\\
&=\ln C+\int_0^{\infty}\left[\ln(u+C)-\ln u\right]\nu(du,z)
\notag\\
&=\ln C+\int_0^{\infty}\ln(1+\frac{u}C)\nu(du,z)
-\int_0^{\infty}\ln u \nu(du,z)
\end{align}
After differentiation we get
\begin{equation}\label{rel04}
\frac{\partial }{\partial s}\int_0^{\infty}\ln u\nu(du,z)
=\frac{\partial }{\partial s}\int_0^{\infty}\ln(1+\frac uC)\nu(du,z)-
\int_{-C}^0\frac{\partial }{\partial s}y(x)dx.
\end{equation}
Relations (\ref{rel01})--(\ref{rel04}) together imply the result.
\end{proof}

\section{The smallest singular value}\label{singular}
In this Section we prove a bound for the minimal singular value of the matrices $\bold X-z\bold I$.
 A corresponding bound for sparse matrices we shall give in the Appendix.
Let $\bold X=\frac1{\sqrt
n}\left(X_{jk}\right)_{j,k=1}^n$ be an $n\times n$ matrix
with i.i.d. entries $X_{jk}$, $j,k=1,\ldots,n$ and $\varepsilon_{jk}$
$j,k=1\ldots,n$ Bernoulli i.\ i.\ d. random variables independent on
$X_{jk}$, $j,k=1,\ldots,n$ with $p_n=\Pr\{\varepsilon_{jk}=1\}$. Assume
that $\E X_{jk}=0$ and $\E X_{jk}^2=1$. We prove the following result.
Denote by $s_1(z)\ge\ldots\ge s_n(z)$ the
singular values of the matrix $\bold X(z):=\bold
X-z\bold I$.
\begin{thm}\label{thm1}Let $X_{jk}$ be independent random variables with sub-Gaussian tails, i. e.
\begin{equation}
\Pr\{|X_{jk}|>t\}\le \exp\{-ct^2\}.
\end{equation}
Then for any $z\in\mathbb C$ such that $|z|\le 4$ and for any
$\gamma>\frac c{\sqrt n}$
\begin{align}
\Pr\{s_n\le\gamma/Cn^{2}\}\le \gamma,
\end{align}
for some positive constants $C$ and $c$.
\end{thm}
The proof of this theorem is  based on the arguments of Rudelson \cite{rud:06}. He
proved the same result for $z=0$ and for a real matrix $\bold X$. To generalize this result  to complex
 $z$
and complex matrices
we need some modifications of his proof. To bound the smallest singular
value in our case we need to consider the complex unit sphere $\mathcal
S^{(n-1)}$ in $\mathbb C^n$. 

\mc{By the symbols $C$ and $c$ with or without indices or without it we shall denote some absolute constants. We shal adapt  Rudelson's
enumeration of constants, i. e. the lower indices of constants correspond the number of the Theorems in 
Rudelson's paper.}

Let $\bold\alpha=(\alpha_1,\ldots,\alpha_n)$ denote a vector in $\mathcal S^{(n-1)}$ in $\mathbb C^n$.
Then $\bold a =(|\alpha_1|,\ldots,|\alpha_n|)$ is an element of the unit sphere $S^{(n-1)}\subset\mathbb R^n$.
We shall use the arguments of Rudelson for real vectors $\bold a$.
Furthermore, we need some modifications of his concentration results for complex random variables.
These are Theorem 3.5  and Lemma 4.2 in \cite{rud:06}.
We start with Theorem 3.5. We may reformulate it as follows.
\begin{thm}\label{rudel}
Let $\beta$ a complex random variable such that $\E\beta=0$ and
$\Pr\{|\beta|>c\}\ge c'$, for some $c,c'>0$. Let $\beta_1,\ldots,\beta_n$
be independent copies of $\beta$. Let $\Delta>0$ and let $\bold
x=(x_1,\ldots,x_n)\in\mathbb C^n$ be a vector such $a<|x_j|< \overline
C_{3.5}a$ for a some $a>0$ and for some positive constant $\overline C_{3.5}$. 
Let  $\varepsilon_j$ be i. i. d. Bernoulli
random variables independent on $\beta_j$, $j=1,\ldots,n$.
Then  there exists a constant $C_{3.5}$ such that for any
$\Delta<\frac a{2\pi}$, for any $j_0=1,\ldots, n$ and any $u, v\in
\mathbb C$
\begin{equation}
\Pr\left\{\left|\sum_{j=1}^n\varepsilon_j\beta_jx_j-\sqrt {np_n}
ux_{j_0}-\sqrt {np_n}v\right|<\Delta\right\}\le
\frac{C_{3.5}}{(np_n)^{\frac52}}\sum_{k=1}^{\infty}P_k^2(x,\Delta),
\end{equation}
where
\begin{equation}\notag
P_k(x,\Delta)=\#\{j:\ |x_j|\in(k\Delta,(k+1)\Delta]\}.
\end{equation}
\end{thm}
\begin{proof}The proof of this Theorem is based on Lemma 3.1 in \cite{rud:06}.
We reformulate this result for the complex  case
\begin{lem}\label{lemrudel}Let $c>0$, $0<\Delta<\frac a{2\pi}$, and $\beta_1,\ldots,\beta_n$ be independent
complex random variables such that $\E\beta_j=0$ and
$\Pr\{|\widetilde\beta_j|>\sqrt 2a\}\ge c$, where $\widetilde
\beta_j=\beta_j-\beta_j'$ and $\beta_j'$ is an  independent copy of
$\beta_j$. Let  $\varepsilon_j$ be i. i. d. Bernoulli random variables
independent on $\beta_j$, $j=1,\ldots,n$. Then, tehere exist constants $c, c'$ such that 
for any $v\in\mathbb C$,
\begin{equation}
\Pr\left\{\left|\sum_{j=1}^n\varepsilon_j\beta_j-v\right|<\Delta\right\}\le
\frac{C_{3.5}}{(np_n)^{\frac52}}\int_{\frac{3a}2}S_{\Delta}^2(y)dy
+ce^{-c'np_n},
\end{equation}
where
\begin{equation}
S_{\Delta}(y)=\sum_{j=1}^n\Pr\{|\widetilde \beta_j|\in[y-\pi\Delta,y+\pi\Delta]\}.
\end{equation}
\end{lem}
\begin{proof}
Let $\beta_j=\xi_{j}+i\ \eta_{j}$, and $v=c+i\ d$.
In this notation we have
\begin{align}
\Pr\left\{\left|\sum_{j=1}^n\beta_j-v\right|<\Delta\right\}
\le
\min\left\{  \Pr\left\{\left|\sum_{j=1}^n\xi_j
-c\right|<\Delta\right\},\right.\notag
\left.\Pr\left\{\left|\sum_{j=1}^n\eta_j
-d\right|<\Delta\right\}\right\}=:Q
\end{align}
Note that
\begin{equation}
|\xi_j|^2+|\eta_j|^2=|\beta_j|^2,
\end{equation}
implies
\begin{equation}
\max\{|\xi_j|,\ |\eta_j|\}\ge\frac{|\beta_j|}{\sqrt 2}.
\end{equation}
By  the Lemma of Ess\'een (see, for example, \cite{Petrov} Lemma 3, p. 38),
for any $v\in \mathbb C$  we have
\begin{equation}
Q\le
C\min\left\{\int_{[-\frac{\pi}2,\frac{\pi}2]}|\phi^{\varepsilon}(t/\Delta)|dt,
\int_{[-\frac{\pi}2,\frac{\pi}2]}|\psi^{\varepsilon}(t/\Delta)|dt\right\},
\end{equation}
where
$\phi^{\varepsilon}(t):=\E\exp\{it\sum_{j=1}^n{\varepsilon_j}\xi_j\}$ and
$\psi^{\varepsilon}(t):=\E\exp\{it\sum_{j=1}^n{\varepsilon_j}\eta_j\}$.
Let $\widetilde\xi_j=\xi_j-(\xi_j)'$ and
$\widetilde\eta_j=\eta_j-(\eta_j)'$ where $(\xi_j)'$ and $(\eta_j)'$
denote independent copies of $\xi_j$ and $\eta_j$ respectively. Note that
\begin{equation}
|\widetilde\xi_j|^2+|\widetilde\eta_j|^2=|\widetilde\beta_j|^2
\end{equation}
This implies that at least $\frac n2$ of the random variables $\widetilde\xi_j$ or $\widetilde \eta_j$,
$j=1,\ldots,n$, satisfy the inequality
\begin{equation}
|\widetilde\xi_j|\ge \frac1{\sqrt 2}|\widetilde \beta_j|,\text{  or  }
|\widetilde \eta_j|\ge \frac1{\sqrt 2}|\widetilde \beta_j|.
\end{equation}
Without loss of generality we shall assume that $m\ge\left[\frac n2\right]$ random variables
$\widetilde\xi_j$ satisfy the
inequality
\begin{equation}
|\widetilde\xi_j|\ge \frac1{\sqrt 2}|\widetilde \beta_j|.
\end{equation}
The last inequality yields
\begin{equation}\label{abc}
\Pr\{|\widetilde\xi_j|\ge a\}\ge \overline c>0.
\end{equation}
Following Rudelson, we introduce the random variable $\tau_j$ by
conditioning on $|\widetilde\xi_j|>2a$. We may repeat from here on his
proof of Lemma 3.1 and Theorem 4.1 in \cite{rud:06} to obtain the result
of Theorem \ref{rudel}.  After simple calculations we get
\begin{equation}\label{new}
|\phi^{\varepsilon}(t)|\le \exp\{-p_n(1-p_n)\sum^*(1-\re\phi_j(t))
-\frac12p_n^2\sum^*(1-|\phi_j(t)|^2)\},
\end{equation}
where $\sum\limits^*$ denote the summation over all indexes $j=1,\ldots,
n$ such that inequality (\ref{abc}) holds and $\phi_j(t)=\E
\exp\{it\xi_j\}$. Furthermore, for all $j$ such that (\ref{abc}) holds we
have
\begin{equation}\label{new1}
1-|\phi_j(t)|^2\ge\overline c\ \E(1-\cos \tau_jt).
\end{equation}
Inequalities (\ref{new}) and (\ref{new1}) together imply
\begin{equation}
|\phi(t)|\le \exp\{-c'f(t)\},
\end{equation}
where
\begin{equation}\notag
f(t)=\E\sum^*(1-\cos\tau_jt).
\end{equation}
In the what follows we repeat  Rudelson arguments for the rest of proof. Let
\begin{align}
T(l,r)&=\{t:\ f(t/\Delta))\le l, |t|\le r\},\\
M&=\max_{|t|\le \frac{\pi}2}f(t/\Delta).
\end{align}
To estimate $M$ from below, notice that
\begin{align}
M&=\max_{|t|\le \frac{\pi}2}f(t/\Delta)\ge \frac1{\pi}
\int_{-\frac{\pi}2}^{\frac{\pi}2}
\E\sum^*(1-\cos(\tau_j/\Delta)t)dt\\&=
\E\sum^*\left(1-\frac2{\pi}\frac{\sin(\tau_j/\Delta)\pi/2}
{\tau_j/\Delta}\right)\ge cm\ge c'n,
\end{align}
since $|\tau_j|/\Delta>2a/\Delta>4\pi$.
We shall use the following result from  Rudelson \cite{rud:06}.
\begin{lem}\label{rud3.2}Let $0<l<M/4$. Then
\begin{equation}
|T(l,\pi/2)|\le c\sqrt{\frac lM}|T(M/4,\pi)|.
\end{equation}
\end{lem}
We have
\begin{align}
Q&\le C\int_{[-\pi/2,\pi/2]}|\phi(t/\Delta)|dt\le C\int_{[-\pi/2,\pi/2]}
\exp\{-c'f(t)\}\notag\\
&\le C\int_0^n|T(l,\pi/2)|e^{-c'l}dl.
\end{align}
 According to the last lemma we get
\begin{equation}
Q\le \frac {C'}{\sqrt M}|T(\frac M4,\pi)|+ce^{-\frac {C'M}{16}}\le
\frac{C'}{\sqrt M}|T(\frac M4,\pi)|+e^{-c'm}.
\end{equation}
Repeating the arguments of Rudelson in \cite{rud:06}, we obtain
\begin{equation}
Q\le \frac C{n^{\frac52}\Delta}\int_{\mathbb R\setminus[-3a/2,3a/2]}
\left(\sum^*\Pr\left\{\tau_j\in
[z-\pi\Delta,z+\pi\Delta]\right\}\right)^2dz+ce^{-c'n}.
\end{equation}
Since $\tau_j$ are symmetric we may change the interval of
 integration set in the previous inequality to $(3a/2,\infty)$.
Moreover, if $z\in(3a/2,\infty)$
\begin{equation}
\Pr\{\tau_j\in[z-\pi\Delta,z+\pi\Delta]\}\le \frac1c
\Pr\{\widetilde \xi_j\in[z-\pi\Delta,z+\pi\Delta]\}\le
\frac1c\Pr\{|\widetilde \beta_j|\in[z-\pi\Delta,z+\pi\Delta]\}.
\end{equation}
Furthermore,
\begin{equation}
1-\re\{\phi_j(t)\}=\E(1-\cos\{\xi_jt\}.
\end{equation}
 repeating the previous arguments, we conclude the proof of Lemma
\ref{lemrudel}.
\end{proof}
We continue to prove Theorem \ref{rudel}.
Recall that $x_j=a_j+i\ b_j$, $u=c+i\ d$ and $v=f+i\ g$. Then the following inequality holds
\begin{align}
&\Pr\{|\sum_{j=1}^n\beta_jx_j-\sqrt nux_{j_0}-\sqrt nv|\le\Delta\}\notag\\&
\qquad\qquad\qquad\le \min\left\{\Pr\{|\sum_{j=1}^n(\xi_ja_j-\eta_jb_j)
-\sqrt n (c\xi_{j_0}-d\eta_{j_0})-\sqrt nf|\le\Delta\},\right.\notag
\\&\qquad\qquad\qquad\qquad\qquad
\left.\Pr\{|\sum_{j=1}^n(\eta_ja_j+\xi_jb_j)
-\sqrt n (c\eta_{j_0}+d\xi_{j_0})-\sqrt nd|\le\Delta\}\right\}.
\end{align}
Note that
\begin{equation}
|\eta_ja_j+\xi_jb_j|^2+|\xi_ja_j-\eta_jb_j|^2=|x_j|^2|\beta_j|^2,
\end{equation}
implies again
\begin{equation}
\max\{|\eta_ja_j+\xi_jb_j|,|\xi_ja_j-\eta_jb_j|\}\ge|x_j||\beta_j|/\sqrt 2.
\end{equation}
Conditioning given $\beta_{j_0}$, we may apply the result of Lemma \ref{lemrudel}.
We obtain
\begin{equation}
\Pr\{|\sum_{j=1}^n\beta_jx_j-\sqrt nux_{j_0}-\sqrt nv|\le\Delta\}\le \frac {C}{m^{\frac52}\Delta}
F(\mu)+ce^{-c'm},
\end{equation}
where
\begin{align}
F(\mu)&=\int_{\frac{3a}2}^{\infty}\widetilde S_{\Delta}^2(y)dy,\notag\\
\widetilde S_{\Delta}(y)&=\sum_{j=1}^n\mu\left(\frac1{|x_j|}[y-\pi\Delta,y+\pi\Delta]\right),
\end{align}
and $\mu$ denotes the distribution of $|\widetilde \beta|$.
Since
\begin{equation}
F(\mu)\le C\Delta\sum_{k=1}^{\infty}|\{j:\ |x_j|\in(k\Delta,(k+1)\Delta]\}|^2,
\end{equation}
we obtain
\begin{equation}
\Pr\{|\sum_{j=1}^n\beta_jx_j-\sqrt nux_{j_0}-\sqrt nv|\le\Delta\}\le \frac {C}{m^{\frac52}}
\sum_{k=1}^{\infty}|\{j:\ |x_j|\in(k\Delta,(k+1)\Delta]\}|^2.
\end{equation}
This completes the proof.
\end{proof}
We also  need  the following lemma.
\begin{lem}Let $x_j=a_j+i\ b_j$, $v=c+i\ d$, $\beta_j=\xi_j+i\ \eta_j$.
Let $\beta$ be a random variable such that $\E\beta=0$, $\E|\beta|^2=1$ and let $\beta_1,\ldots,\beta_n$,
be independent copies of $\beta$. Let $0<r<R$ and let $x_1,\ldots,x_m\in\mathbb C$ such that
$\frac r{\sqrt m}\le |x_j|\le \frac R{\sqrt m}$ for any $j$.
Then there exist  constants $ C_{4.2}$ and $c_{4.2}$ such that 
for any $t>\frac{c_{4.2}}{\sqrt m}$ and for any $v\in\mathbb C$
\begin{equation}
\Pr\left\{\left|\sum_{j=1}^n\beta_jx_j-v\right|<t\right\}\le C_{4.2}t
\end{equation}
\end{lem}
\begin{proof}
We use the simple inequality
 \begin{equation}
\Pr\left\{\left|\sum_{j=1}^n\beta_jx_j-v\right|<t\right\}\le\min\{A,B\},
\end{equation}
where
\begin{align}\notag
A&=
\Pr\left\{\left|\sum_{j=1}^n(\xi_ja_j-\eta_jb_j)-c\right|<t
\right\}\\
B=&\Pr\left\{\left|\sum_{j=1}^n(\eta_ja_j+\xi_jb_j)-d\right|<t\right\}
\end{align}
Note that random variables $\widehat \xi_j=\xi_ja_j-\eta_jb_j$ (resp. $\overline \xi_j=\xi_ja_j-\eta_jb_j$)
are independent for $j=1,\ldots, n$,
\begin{equation}
\max\{\sum_{j=1}^m\E|\widehat \xi_j|^3,\ \sum_{j=1}^m\E|\overline \xi_j|^3\}\le \frac{CR^3}{\sqrt m},
\end{equation}
and
\begin{equation}
\max\{\sum_{j=1}^m\E||\widehat \xi_j|^2,\ \sum_{j=1}^m\E|\overline \xi_j|^2\}
\ge \frac{r^2\sigma^2}{\sqrt 2}.
\end{equation}
Applying the Berry--Ess\'een inequality, we obtain the result.
\end{proof}

To conclude the proof of Theorem \ref{thm1} we repeat the proof of Rudelson \cite{rud:06} in the rest.

\section{Proof of the main Theorem}\label{proof}
In this Section we give the proof of Theorem \ref{thm0}. The proof of Theorem \ref{sparse} is similar.
We  have to use Theorem \ref{thm1.1} instead of Theorem \ref{thm1} and instead of Bai's results 
we may use the result of Section \ref{convergence} for $z=0$ only. For any $z\in \mathbb C$ we introduce the set
$\Omega_n(z)=\{\omega\in\Omega:\ n^{-3}\ge s_n(\bold X-z\bold I), \ s_1(\bold X)\le 4\}$.
From  Bai \cite{Bai:93b} it follows
that
\begin{equation}
\Pr\{s_1(X)\ge 4\}\le Cn^{-\frac18}.
\end{equation}
 According to Theorem \ref{thm1},
\begin{equation}
\Pr\{n^{-3}\ge s_n(\bold X-z\bold I)\}\le Cn^{-\frac12}.
\end{equation}
These  inequalities imply
\begin{equation}\label{truncation}
\Pr\{\Omega_n(z)\}\le Cn^{-\frac18}.
\end{equation}
Let $r=r(n)$ such that $r(n)\to 0$ as $n\to\infty$.  A more specific choice will be made  later.
Consider  the potential $U_{\mu_n}^{(r)}$.
We have
\begin{align}
U_{\mu_n}^{(r)}&=-\frac1n\E\log|\det(\bold X-z\bold I-r\xi\bold I)|\\&
=-\frac1n\sum_{j=1}^n\E\log|\lambda_j-r\xi-z|I_{\Omega_n(z)}
-\frac1n\sum_{j=1}^n\E\log|\lambda_j-r\xi-z|I_{\Omega_n^{(c)}(z)}\\&=
\overline U_{\mu_n}^{(r)}+\widehat U_{\mu_n}^{(r)},
\end{align}
where
$I_A$ denotes an indicator function of an event $A$ and $\Omega_n^{(c)}(z)$ denotes the 
complement of $\Omega_n(z)$.
\begin{lem}\label{lem5.1}Assuming the conditions of Theorem \ref{thm1}, for $r$ such that $-n^{-1/12}\log r\to 0$ as $n\to\infty$, we have
\begin{equation}\label{0*}
\widehat U_{\mu_n}^{(r)}\to 0,\text{  as  }n\to\infty.
\end{equation}
\end{lem}
\begin{proof}
By definition, we have
\begin{align}\label{1*}
\widehat U_{\mu_n}^{(r)}=-\frac1n\sum_{j=1}^n\E\log|\lambda_j-r\xi-z|I_{\Omega_n^{(c)}(z)}.
\end{align}
Applying Cauchy's inequality, we get, for any $\alpha>0$,
\begin{align}\label{2*}
|\widehat U_{\mu_n}^{(r)}|&\le \frac1n\sum_{j=1}^n\E^{\frac1{1+\alpha}}|\log|\lambda_j-r\xi-z||^{1+\alpha}
\left(\Pr\{\Omega_n\}\right)^{\frac{\alpha}{1+\alpha}}\notag\\&\le
\left(\frac1n\sum_{j=1}^n\E|\log|\lambda_j-r\xi-z||^{1+\alpha}\right)^{\frac1{1+\alpha}}\left(\Pr\{\Omega_n\}\right)^{\frac{\alpha}{1+\alpha}}.
\end{align}
Furthermore, since $\xi$ is uniformly distributed in the unit disc and independent of $\lambda_j$, we may write
\begin{equation}
\E|\log|\lambda_j-r\xi-z||^{1+\alpha}=\frac1{2\pi}\E\int_{|\zeta|\le1}|\log|\lambda_j-r\zeta-z||^{1+\alpha}d\zeta=
\ E J_1^{(j)}+\E J_2^{(j)}+\E J_3^{(j)},
\end{equation}
where
\begin{align}
J_1^{(j)}&=\frac1{2\pi}\int_{|\zeta|\le 1,\ |\lambda_j-r\zeta-z|\le\varepsilon}|\log|\lambda_j-r\zeta-z||^{1+\alpha}d\zeta\\
J_2^{(j)}&=\frac1{2\pi}\int_{|\zeta|\le 1,\ \frac1{\varepsilon}>|\lambda_j-r\zeta-z|>\varepsilon}|\log|\lambda_j-r\zeta-z||^{1+\alpha}d\zeta\\
J_3^{(j)}&=\frac1{2\pi}\int_{|\zeta|\le 1,\ |\lambda_j-r\zeta-z|>\frac1{\varepsilon}}|\log|\lambda_j-r\zeta-z||^{1+\alpha}d\zeta
\end{align}
Note that
\begin{equation}
|J_2^{(j)}|\le\log\left(\frac1{\varepsilon}\right).
\end{equation}
Since for any $b>0$, the function $-u^a\log u$ is not decreasing on the interval $[0,\exp\{-\frac1{b}]$,
we have for $0<u\le\varepsilon<\exp\{-\frac1{b}\}$,
\begin{equation}
-\log u\le \varepsilon^{b}u^{-b}\log\left(\frac1{\varepsilon}\right).
\end{equation}
Using this inequality, we obtain, for $b(1+\alpha)<2$,
\begin{align}\label{finish1}
|J_1^{(j)}|&\le\frac1{2\pi}\varepsilon^{b(1+\alpha)}\left(\log\left(\frac1{\varepsilon}\right)\right)^{1+\alpha}
\int_{|\zeta|\le 1,\ |\lambda_j-r\zeta-z|\le\varepsilon}|\lambda_j-r\zeta-z|^{-b(1+\alpha)}d\zeta\\
&\le\frac1{2\pi r^2}\varepsilon^{b}\log\left(\frac1{\varepsilon}\right)\int_{|\zeta|\le \varepsilon}|\zeta|^{-b(1+\alpha)}d\zeta\le C(\alpha, b)\varepsilon^{2}r^{-2}\left(\log\left(\frac1{\varepsilon}\right)\right)^{1+\alpha}
\end{align}
If we choose $\varepsilon=r$, then we get
\begin{equation}\label{finish2}
|J_1^{(j)}|\le C(\alpha, b)\left(\log\left(\frac1{r}\right)\right)^{1+\alpha}
\end{equation}
The following bound holds for $\frac1n\sum_{j=1}^n\E J_3^{(j)}$. Note that
$|\log x|^{1+\alpha}\le \varepsilon^2|\log\varepsilon|^{1+\alpha}x^2$ for $x\ge\frac1{\varepsilon}$ and 
sufficiently small $\varepsilon$.
Using this inequality, we obtain
\begin{align}\label{finish3}
\frac1n\sum_{j=1}^n\E J_3^{(j)}\le C(\alpha)\varepsilon^2|\log\varepsilon|\frac1n\sum_{j=1}^n \E |\lambda_j-r\zeta-z|^2
\le C(\alpha)(1+|z|^2+r^2)\varepsilon^2|\log\varepsilon|\notag\\
\le C(\alpha)(2+|z|^2)r^2|\log r|.
\end{align}

The inequalities (\ref{finish1})--(\ref{finish3}) together imply that
\begin{equation}\label{3*}
|\frac1n\sum_{j=1}^n\E|\log|\lambda_j-r\xi-z||^{1+\alpha}|\le
C\left(\log\left(\frac1{r}\right)\right)^{1+\alpha}.
\end{equation}
 Furthermore, the inequalities (\ref{1*}), (\ref{2*}), and (\ref{3*}) together imply
\begin{equation}
|\widehat U_{\mu_n}^{(r)}|\le C\left(\log\left(\frac1{r}\right)\right)n^{-\frac{\alpha}{3(1+\alpha)}}
\end{equation}
We choose $\alpha=3$  and rewrite the last inequality as follows
\begin{equation}
|\widehat U_{\mu_n}^{(r)}|\le C\left(\log\left(\frac1{r}\right)\right)n^{-\frac{1}{4}}
\end{equation}
If we choose $r$ such that $\log(1/r)n^{-1/4}\to 0$, then (\ref{0*}) holds. Thus the 
Lemma is proved.
\end{proof}
We shall investigate $\overline U_{\mu_n}^{(r)}$ now.
We may write
\begin{align}
\overline U_{\mu_n}^{(r)}&=-\frac1n\sum_{j=1}^n\E\log|\lambda_j-z-r\xi|I_{\Omega_n}=
-\frac1n\sum_{j=1}^n\E\log(s_j(\bold X(z,r))I_{\Omega_n}\\
&=-\int_{n^{-3}}^{4+|z|}\log xd\E\overline F_n(x,z,r),
\end{align}
where $\overline F_n(\cdot,z,r)$ is  the distribution function corresponding  to the 
restriction of the measure $\nu_n(\cdot,z,r)$ on the set $\Omega_n$.
Introduce the notation
\begin{equation}
\overline U_{\mu}=-\int_{n^{-3}}^{4+|z|}\log x dF(x,z)
\end{equation}
Integrating  by parts, we get
\begin{equation}
\overline U_{\mu_n}^{(r)}-\overline U_{\mu}=-\int_{n^{-3}}^{4+|z|}\frac{\E F_n(x,z,r)-F(z,r)}x dx+
C\theta\sup_x|\E F_n(x,z,r)-F(z,r)|\log n,
\end{equation}
where $\theta$ denotes some constant such that $|\theta|\le 1$.
This implies that
\begin{equation}\label{009}
|\overline U_{\mu_n}^{(r)}-\overline U_{\mu}|\le C\log n\sup_x|\E F_n(x,z,r)-F(x,z)|.
\end{equation}
Note that, for any $r>0$, $|s_j(z)-s_j(z,r)|\le r$. This implies that
\begin{equation}
\E F_n(x-r,z)\le \E F_n(x,z,r)\le \E F_n(x+r,z).
\end{equation}
 Hence, we get
\begin{equation}
\sup_x|\E F_n(x,z,r)-F(x,z)|\le \sup_x|\E F_n(x,z)-F(x,z)|+\sup_x|F(x+r,z)-F(x,z)|.
\end{equation}
Since  the distribution function $F(x,z)$ has a density $p(x,z)$  which is bounded
(see Remark \ref{rem:01})
we obtain
\begin{equation}\label{08}
\sup_x|\E F_n(x,z,r)-F(x,z)|\le \sup_x|\E F_n(x,z)-F(x,z)|+Cr.
\end{equation}
Choose $r=cn^{-\frac18}$.
Inequalities \ref{08} and (\ref{supremum}) together imply
\begin{equation}\label{09}
\sup_x|\E\overline F_n(x,z,r)-\overline F(x,z)|\le Cn^{-\frac1{8}}.
\end{equation}
From inequalities (\ref{09}) and (\ref{009}) it follows that
\begin{equation}
|\overline U_{\mu_n}^{(r)}-\overline U_{\mu}|\le Cn^{-\frac1{8}}\log n.
\end{equation}
Note that
\begin{equation}
|\overline U_{\mu_n}^{(r)}-U_{\mu}|\le |\int_0^{n^{-3}}\log xdF(x,z)|\le Cn^{-\frac18}\log n.
\end{equation}

Let $K=\{z\in\mathbb C :\ |z|\le 4\}$ and let $K^{(c)}$ denote  $\mathbb C\setminus K$.
According to inequality (\ref{thm0}), we have
\begin{equation}\label{rep1}
1-p_n:=\E\mu_n^{(r)}(K^{(c)})\le\Pr\{s_1(\bold X)>3\}\le \sup_x|F_n(x)-M_1(x)|\le Cn^{-\frac18}.
\end{equation}
Furthermore, let ${\overline{\mu}}_n^{(r)}$ and ${\widehat{\mu}}_n^{(r)}$ be  probability measures supported on the compact set
$K$ and $K^{(c)}$ respectively, such that
\begin{equation}\label{rep2}
\E\mu_n^{(r)}=p_n{\overline{\mu}}_n^{(r)}+(1-p_n){\widehat{\mu}}_n^{(r)}.
\end{equation}
 Introduce the logarithmic potential of  the measure ${\overline{\mu}}_n^{(r)}$,
\begin{equation}
U_{{\overline{\mu}}_n^{(r)}}=-\int\log|z-\zeta|d{\mu_n{\overline{\mu}}_n^{(r)}(\zeta)}.
\end{equation}
Similar to  the proof of Lemma \ref{lem5.1} we show that
\begin{equation}
\lim_{n\to \infty}|U_{\mu_n}^{(r)}-U_{{\overline{\mu}}_n^{(r)}}|\le Cn^{-\frac18}\log n.
\end{equation}
This implies that
\begin{equation}
\lim_{n\to\infty}U_{{\overline{\mu}}_n^{(r)}}(z)=U_{\mu}(z)
\end{equation}
for all $z\in\mathbb C$.
Since the measures ${\overline{\mu}}_n^{(r)}$ are compactly supported, Theorem 6.9 from \cite{saff}
and Corollary 2.2 from \cite{saff}
 (see also the Appendix, Theorem \ref{LEP} and Corollary \ref{cor2.2}), together  imply that
\begin{equation}\label{rep3}
\lim_{n\to \infty}\overline{\mu}_n^{(r)}=\mu
\end{equation}
in the weak topology.
 Inequality (\ref{rep1}) and relations (\ref{rep2}) and (\ref{rep2}) together imply
that
\begin{equation}
 \lim_{n\to\infty}\E\mu_n^{(r)}=\mu
\end{equation}
in weak topology.
 Finally, by Lemma \ref{smo1} we get
\begin{equation}
 \lim_{n\to\infty}\E\mu_n=\mu
\end{equation}
in the weak topology.
Thus Theorem \ref{thm0} is proved.
\section{Appendix}In this Section we  collect some technical results.
\begin{lem}\label{ap}Let $\varkappa_3=\max_{j,k}\E|X_{jk}|^3$. The following inequality holds
\begin{align}
\frac1{n\sqrt n}\sum_{j,k=1}^n\E|X_{jk}|(|T^{(jk)}_{k+n,j}|+|T^{(jk)}_{j,k+n}|)\le
\frac{C\varkappa_3}{\sqrt nv^3}
\end{align}
\end{lem}
\begin{proof}

Introduce  the notations
\begin{equation}
B:=\frac1{n\sqrt n}\sum_{j,k=1}^n\E|X_{jk}|(|T^{(jk)}_{k+n,j}|+|T^{(jk)}_{j,k+n}|)
\end{equation}
and
\begin{align}
B_1&:=\frac1{n^2\sqrt n}\sum_{j,k=1}^n\E|X_{jk}|^3|R^{(jk)}_{k+n,j}|^2|R_{k+n,j}|\notag\\
B_2&:=\frac1{n^2\sqrt n}\sum_{j,k=1}^n\E|X_{jk}|^3|R^{(jk)}_{k+n,j}||R^{(jk)}_{k+n,k+n}||R_{j,j}|\notag\\
B_3&:=\frac1{n^2\sqrt n}\sum_{j,k=1}^n\E|X_{jk}|^3|R^{(jk)}_{k+n,k+n}||R^{(jk)}_{j,j}||R_{k+n,j}|\notag\\
B_4&:=\frac1{n^2\sqrt n}\sum_{j,k=1}^n\E|X_{jk}|^3|R^{(jk)}_{k+n,k+n}||R^{(jk)}_{k+n,j}||R_{j,j}|\notag\\
B_5&:=\frac1{n^2\sqrt n}\sum_{j,k=1}^n\E|X_{jk}|^3|R^{(jk)}_{j,k+n}||R^{(jk)}_{k+n,j}||R_{j,k+n}|\notag\\
B_6&:=\frac1{n^2\sqrt n}\sum_{j,k=1}^n\E|X_{jk}|^3|R^{(jk)}_{j,j}||R^{(jk)}_{k+n,k+n}||R_{j,k+n}|\notag\\
B_7&:=\frac1{n^2\sqrt n}\sum_{j,k=1}^n\E|X_{jk}|^3|R^{(jk)}_{j,k+n}||R^{(jk)}_{j,j}||R_{k+n,j}|\notag\\
B_8&:=\frac1{n^2\sqrt n}\sum_{j,k=1}^n\E|X_{jk}|^3|R^{(jk)}_{j,k+n}|^2|R_{j,k+n}|\notag\\
\end{align}
It is easy to check that
\begin{equation}
\max\{B_k,\ k=1,\ldots,8\}\le \frac{C\varkappa_3}{\sqrt nv^3}.
\end{equation}
This implies that
\begin{equation}
B\le\frac{C\varkappa_3}{\sqrt nv^3}.
\end{equation}

\end{proof}

\begin{lem}\label{ap1}
Let $\mu_n$ be the empirical spectral measure of  the matrix $\bold X$ and
$\nu_r$ be the uniform distribution on the disc of radius $r$. Let
$\mu_n^{(r)}$ be  the empirical spectral measure of the matrix $\bold X(r)= \bold
X-r\xi\bold I$, where $\xi$ is a random variable which is  uniformly
distributed  on the unit disc. Then the measure $\E\mu_n^{(r)}$ is the
convolution of the measures $\E\mu_n$ and $\nu_r$, i. e.
\begin{equation}
\E\mu_n^{(r)}=(\E\mu_n)*(\nu_r).
\end{equation}
\end{lem}
\begin{proof}
Let $J$ be a random variable which is  uniformly distributed on the set
$\{1,\ldots,n\}$. Let $\lambda_1,\ldots,\lambda_n$ be  the eigenvalues of the
matrix $\bold X$. Then $\lambda_1+r\xi,\ldots,\lambda_n+r\xi$ are
eigenvalues of the matrix $\bold X(r)$. Let $\delta_x$ be denote the Dirac
measure. Then
\begin{equation}
\mu_n=\frac1n\sum_{j=1}^n\delta_{\lambda_j}
\end{equation}
and
\begin{equation}
 \mu_n^{(r)}=\frac1n\sum_{j=1}^n\delta_{\lambda_j+r\xi}.
 \end{equation}
 Denote by $\mu_{nj}$  the distribution of $\lambda_j$.
 Then
 \begin{equation}
 \E\mu_n=\frac1n\sum_{j=1}^n\mu_{nj}
 \end{equation}
 and
 \begin{equation}
 \E\mu_n^{r}=\frac1n\sum_{j=1}^n\mu_{nj}*\nu_r=
 \left(\frac1n\sum_{j=1}^n\mu_{nj}\right)*(\nu_r)=(\E\mu_n)*(\nu_r).
 \end{equation}
  The Lemma is proved.
\end{proof}
 Let
\begin{equation}
f_n^{(r)}(t,v)=\int_{-\infty}^{\infty}\int_{-\infty}^{\infty}\exp\{itx+
ivy\} dG_n^{(r)}(x,y)
\end{equation}
and
\begin{equation}
f_n(t,v)=\int_{-\infty}^{\infty}\int_{-\infty}^{\infty}\exp\{itx+
ivy\}dG_n(x,y),
\end{equation}
where
\begin{equation}
G_n^{(r)}(x,y)=\frac1n\sum_{j=1}^n\Pr\{\re{\lambda_j+r\xi}\le
x,\im{\lambda_j+r\xi}\le y\},
\end{equation}
and
\begin{equation}
G_n(x,y)=\frac1n\sum_{j=1}^n\Pr\{\re{\lambda_j}\le x,\im{\lambda_j}\le
y\}.
\end{equation}
Denote by $h(t,v)$ the characteristic function of  the joint distribution of
 the  real and imaginary parts of $\xi$,
\begin{equation}
h(t,v)=\int_{-\infty}^{\infty}\int_{-\infty}^{\infty}\exp\{iux+ ivy\}
dG(x,y).
\end{equation}

\begin{lem}\label{ap2}
The following relations  hold
\begin{equation}
f_n^{(r)}(t,v)=f_n(t,v)h(rt,rv).
\end{equation}
If for any $t,v$ there exists $\lim_{n\to\infty}f_n(t,v)$, then
\begin{equation}
\lim_{r\to0}\lim_{n\to\infty}f_n^{(r)}(t,v)=\lim_{n\to\infty}\lim_{r\to0}
f_n^{(r)}(t,v)=\lim_{n\to\infty}f_n(t,v).
\end{equation}
\end{lem}
\begin{proof}The first equality follows immediately from the  independence of the random variable 
$\xi$ and  the
matrix $\bold X$. Since $\lim_{r\to0}h(rt,rv)=h(0,0)=1$ the first
equality implies the second one.
\end{proof}

\begin{lem}\label{ap3} Let $F$ and $G$ be  distribution functions  with
 Stieltjes transforms $S_F(z)$
and $S_G(z)$ respectively. Assume that
$\int_{-\infty}^{\infty}|F(x)-G(x)|dx<\infty$. Let $G(x)$ have a bounded
support $J$ and  density bounded by some constant $K$. Let $V>v_0>0$ and
$a$ be positive numbers such that
$$
\gamma=\frac1{\pi}\int_{|y|\le a}\frac1{u^2+1}\ du > \frac34.
$$
 Then there exist some
constants $C_1,\, C_2,\, C_3$ depending on $J$ and $K$ only such that
\begin{align}
\sup_x|F(x)-G(x)|&\le \, C_1\ \sup_{x\in J} \
\int_{-\infty}^x|S_F(u+iV)-S_G(u+iV)|\ du\notag\\&+ \sup_{u\in J}
\int_{v_0}^V |S_F(u+iv)-S_G(u+iv)|dv+ C_3\ v_0
\end{align}
\end{lem}
\subsection{Some facts from logarithmic potential theory}\label{potential}
We cite here some definitions and Theorems  about  logarithmic  potentials, see \cite{saff}.
Let $\Sigma\subset\mathbb C$ be a compact set of the complex plane and $\mathcal M(\Sigma)$
the collection of all positive  Borel  probability measures with support in $\Sigma$.
The {\it logarithmic energy} of  $\mu\in\mathcal M(\Sigma)$ is defined as
\begin{equation}
I(\mu):=\int\int\log\frac{1}{|z-t|}d\mu(z)d\mu(t),
\end{equation}
and the energy of $\Sigma$ by
\begin{equation}
V:=\inf\{I(\mu)|\mu\in\mathcal M(\Sigma)\}.
\end{equation}
The quantity
\begin{equation}
cap(\Sigma):=e^{-V}
\end{equation}
is called the {\it logarithmic capacity} of $\Sigma$.

The {\it capacity} of an arbitrary  Borel set $E$ is defined as
\begin{equation}
\cap(E):=sup\{cap(K)|K\subset E, K\text{ compact}\}.
\end{equation}
Note that every Borel set of capacity zero has zero two-dimensional Lebesgue measure.
A property is said to hold {\it quasi-everywhere} (q. e.) on a set $E$ if the set of exceptional points
is of capacity zero.
The next Theorem is called {\it Lower Envelope Theorem}
\begin{thm}\label{LEP}
Let ${\mu_n}$, $n=1,2\ldots$, be a sequence of positive  Borel  probability measures  having 
support in a
fixed compact set.
If $\mu_n\to\mu$ weakly, then
\begin{equation}
\liminf_{n\to\infty}U^{\mu_n}(z)=U^{\mu}(z)
\end{equation}
for quasi-every $z\in\mathbb C$.
\end{thm}
The following fact is  Corollary 2.2 from the Unicity Theorem of logarithmic potential theory
(see \cite{saff}, p. 98).
\begin{cor}\label{cor2.2}
If $\mu$ and $\nu$ are compactly supported measures and the potentials $U^{\mu}$ and $U^{\nu}$
coincides almost everywhere with respect to two-dimensional Lebesgue measure, then $\mu=\nu$.
\end{cor}
For reader convenience we give here the statement of Theorem 1.2 from \cite{saff}.
\begin{thm}\label{thm_saff}
Let $\mu$ be a finite positive measure of compact support on the plane.
Then for any $z_0$ and $r>0$ the mean value
\begin{equation}
L(U^{\mu};z_0,r):=\frac1{2\pi}\int_{-\pi}^{\pi}U^{\mu}(z_0+r\exp\{i\theta\}d\theta
\end{equation}
exists as a finite number, and $L(U^{\mu};z_0,r)$ is a non-increasing function of $r$ that is absolutely
continuous on any closed subinterval of $(0,\infty)$.
Furthermore,
\begin{equation}
\lim_{r\to0}L(U^{\mu};z_0,r)=U^{\mu}(z_0).
\end{equation}
\end{thm}
\subsection{Minimal singular  values} of sparse matrices  \label{sparse1}
In this Section we reformulate some statements from the paper  of Rudelson  \cite{rud:06} to 
 adapt his proof  to  sparse matrices. Let $\varepsilon_{jk}$ 
be independent Bernoulli random variables with 
$\Pr\{\varepsilon_{jk}=1\}=p_n$. Assume that $\varepsilon_{jk}$, $j,k=1,\ldots,n$ are independent on 
$X_{jk}$, $j,k=1,\ldots,n$. Consider  the matrix
\begin{equation}
\bold X^{\varepsilon}=(\frac1{\sqrt {p_n}}\varepsilon_{jk}X_{jk})_{j,k=1}^n.
\end{equation}
\begin{thm}\label{thm1.1}Let $X_{jk}$, $j,k=1,\ldots n$  be  centered sub-Gaussian random variables of 
variance 1. 
Then for any $\gamma>c_{1,1}/\sqrt{np_n}$
\begin{equation}
\Pr\left\{\text{there exists } \bold x\in S^{n-1}\ \Big| \ \|\bold X ^{\varepsilon}\bold x\|\le 
\frac{\gamma \sqrt{p_n}}{C_{11}(np_n)^{\frac32}}\right \}\le c\gamma p_n^{-\frac52}
\end{equation}
if $n$ is large enough.
\end{thm}

\mc{The generalization of this result to the complex case is based on similar arguments
as in Section
\ref{singular}  for the case $p_n=1$}.
\begin{proof}
We  adapt Rudelson's proof for sparse matrices
 giving only the neccccessary  new statements of some Lemmas and Theorems in Rudelson's proof. 
To prove   these results is enough to repeat Rudelson's proof of the corresponding Theorems and Lemmas.
\begin{lem}\label{lemma3.1}(Lemma 3.1 in \cite{rud:06})
Let $c>0$, $0<\Delta<a/{2\pi}$ and  let 
$\xi_1,\ldots,\xi_n$ be independent random variables such that 
$\E \xi_i=0$,
$\Pr\{\xi_i>2a\}\ge c$ and $\Pr\{-\xi_i>2a\}\ge c$.  For $y\in\mathbb R$ set
\begin{equation}
S_{\Delta}(y)=\sum_{j=1}^n\frac12(\Pr\{\xi_j\in[y-\pi\Delta,y+\pi\Delta]\}+\Pr\{-\xi_j\in[y-\pi\Delta,y+\pi\Delta]\}),
\end{equation}
 Let $\varepsilon_1,\ldots,\varepsilon_n$ be identically distributed Bernoulli random variables
independent on $\xi_1,\ldots,\xi_n$ and independent 
in aggregate,  with $\Pr\{\varepsilon_j=1\}=p_n$.
Then for any $v\in\mathbb R$
\begin{equation}
\Pr\left\{\Big|\sum\varepsilon_j\xi_j-v\Big|<\Delta\right\}\le
\frac{C}{n^{\frac52}p_n^{\frac32}\Delta}\int_{\frac{3a}2}^{\infty}S_{\Delta}^2(y)dy+c\exp\{-c'np_n\}.
\end{equation}
\end{lem}
\begin{thm}\label{theorem3.5}(Theorem 3.5 in \cite{rud:06})
Let $\xi_1,\ldots,\xi_n$ i.\ i.\ d.\  be sub-Gaussian random variables such that $\E\xi_i=0$ and
$\Pr\{\xi_i>c\}\ge c'$, $\Pr\{-\xi_i>c\}\ge c'$ for some $c,c'>0$. Let $\Delta>0$ and let 
$(x_1,\ldots,x_m)\in\mathbb R^m$ be a vector such $a<|x_j|<\overline C_{3.5}a/\sqrt{p_n}$.
Let $\varepsilon_1,\ldots,\varepsilon_n$ be independent on $\xi_1,\ldots,\xi_n$ and independent 
in aggregate identically distributed Bernoulli random variables with $\Pr\{\varepsilon_j=1\}=p_n$.
Then for any $\Delta<a/(2\pi)$ and for any $v\in \mathbb R$
\begin{equation}
\Pr\left\{\left|\sum_{j=1}^m\xi_j\varepsilon_jx_j-v\right|<\Delta\right\}\le\frac {C_{3.5}}{(mp_n)^{\frac52}}
\sum_{k=1}^{\infty}P_k^2(x,\Delta),
\end{equation}
where
\begin{equation}
P_k^2(x,\Delta)=|\{j\ \Big|\ |x_j|\in(k\Delta,(k+1)\Delta]\}|.
\end{equation}
\end{thm}
\begin{lem}\label{lem4.1}(Lemma 4.1 in \cite{rud:06})
 Assuming the conditions of Theorem \ref{thm1.1}, for the matrix $\bold X^{\varepsilon}$ and for 
every $v\in\mathbb R$,
  we have
\begin{equation}
\Pr\{\|\bold X^{\varepsilon}\|\le C_{4.1}\sqrt{np_n}\}\le\exp\{-c_{4.1}np_n\}.
\end{equation}
\end{lem}
\begin{lem}\label{lem4.2}(Lemma 4.2 in \cite{rud:06})
Let $\xi_1,\ldots,\xi_n$ be i.\ i.\ d.\ sub-Gaussian random variables such that $\E\xi_i=0$ and
$\E\xi_i^2=1$. Let $0<r<R$ and let $x_1,\ldots,x_m\in\mathbb R$ be such that $\frac r{\sqrt m}<|x_j|
<\frac {R}{\sqrt{mp_n}}$ for any $j$. Then for $t\ge \frac{c_{4.2}}{\sqrt {mp_n}}$ and for any $v\in\mathbb R$ 
\begin{equation}
\Pr\left\{\Big|\sum_{j=1}^m\xi_j\varepsilon_jx_j-v\Big|<t\right\}\le C_{4.2}t/\sqrt{p_n}.
\end{equation}
\end{lem}
\begin{lem}\label{lemma4.4}(Lemma 4.4 in \cite{rud:06})Let $\Delta>0$ and let $Y$ be a random variable such that for any $t\ge \Delta$, $\Pr\{|Y-v|<t\}\le Lt$.
Let $\bold y=(Y_1,\ldots,Y_n)$ be a random vector, whose coordinates are independent copies of $Y$.
Then for any $\bold z\in\mathbb R^n$
\begin{equation}
\Pr\left\{\|\bold y-\bold z\|\le\Delta\sqrt n\right\}\le(C_{4.4}L\Delta)^n.
\end{equation}
\end{lem}
We define the set $\sigma(\bold x)$ for any $\bold x\in \mathcal S^{(n-1)}$ as
\begin{equation}
\sigma(\bold x)=\{i\ |\ |x_i|\le R/\sqrt{np_n}\}.
\end{equation}
Let $P_I$ be the coordinate projection on the set $I\subset\{1,\ldots,n\}$.
Set 
\begin{align}
V_P&=\{\bold x\in\mathcal S^{(n-1}\ |\ \|P_{\sigma(\bold x)}\bold x\|<r\}\notag\\
V_S&=\{\bold x\in\mathcal S^{(n-1}\ |\ \|P_{\sigma(\bold x)}\bold x\|\ge r\}\notag.
\end{align}
\begin{lem}\label{lemma5.1}(Lemma 5.1 in \cite{rud:06})For any $r<1/2$
\begin{equation}
\log N(V_P,\mathcal B_2^n,2r)\le \frac{np_n}{R^2}\log{\frac{3R^2}{r{p_n}}}.
\end{equation}
\end{lem}
\begin{lem}\label{lemma5.2}(Lemma 5.2 in \cite{rud:06})
\begin{equation}
\Pr\left\{\text{there exists }\bold x\in V_P\ \Big|\ \|\bold X^{\varepsilon}\bold x\|
\le C_{4.1}\sqrt{np_n}/2\right\}\le \exp\{-c_{4.1}np_n\}.
\end{equation}
\end{lem}
For $\bold x=(x_1,\ldots,x_n)\in V_S$ denote 
\begin{equation}
J(\bold x)=\{j\ |\ \frac r{2\sqrt n}\le |x_j|\le \frac R{\sqrt{np_n}}\}.
\end{equation}
Note that
\begin{equation}
|J(\bold x)|\ge (r^2/2R^2p_n)n=:m.
\end{equation}
Let $0<\Delta<r/2\sqrt n$ be a number to be chosen later.
We shall cover the interval $[\frac r{2\sqrt n},\frac R{\sqrt{np_n}}]$ by
\begin{equation}
k=\left[\frac{R/\sqrt {p_n}-r/2}{\sqrt n\Delta}\right].
\end{equation}
consecutive intervals $(j\Delta,(j+1)\Delta]$, where $j=k_0,(k_0+1),\ldots,(k_0+k)$ and $k_0$ 
is the largest number such that $k_0\Delta<r/2\sqrt n$.

\begin{def}\label{def5.3}(Definition 5.3 in \cite{rud:06}).Let $\Delta>0$ and $Q>1$. We say that a vector $\bold x\in V_S$ has a $(\Delta,Q)$-regular profile if there exists a set $J\subset J(\bold x)$ such that $|J|\ge m/2$ and 
\begin{equation}
\sum_{i=1}^{\infty}P_i^2(\bold x|_J,\Delta)\le Qm^{\frac52}\Delta=:C_{5.3}Q\frac{m^2}k.
\end{equation}
\end{def}
\begin{lem}\label{lem6.1}(Lemma 6.1 in \cite{rud:06})
Let $\Delta\le \frac r{4\pi\sqrt n}$. Let $\bold x\in V_S$ be a vector of $(\Delta,Q)$-regular profile. Then for any $t\ge\Delta$
\begin{equation}
\Pr\left\{\left|\sum_{j=1}^n\xi_j\varepsilon_jx_j-v\right|\le t\right\}\le C_{6.1}Qt/{p_n}.
\end{equation}
\end{lem}
\begin{thm}\label{thm6.2}(Theorem 6.2 in \cite{rud:06})
Let $\frac r{4\pi \sqrt{np_n}}>\Delta>0$ and let $U$ be the set of vectors of $(\Delta,Q)$-regular profile. Then 
\begin{equation}
\Pr\{\text{there exists }\bold x\in U\ |\ \|\bold X^{\varepsilon}\bold x\|\le\frac{\Delta\sqrt{p_n}}
{2\sqrt n}\}\le
C_{6.1}Q\Delta n/\sqrt{p_n}
\end{equation}
\end{thm}
\begin{lem}\label{lem7.1}(Lemma 7.1 in \cite{rud:06})
Let $\overline C_{7.1}(np_n)^{-\frac32}\le\Delta\le {(np_n)^{-\frac12}}$ where $C_{7.1}=\frac{2R^3}{r^2}$
and let $W_S$ be the set of vectors of $(\Delta,Q)$-singular profile.
Let $\eta>0$ be such that 
\begin{equation}
C(\eta)<C_{5.3}Q,
\end{equation}
where $C(\eta)$ is the function defined in Lemma 2.1 in \cite{rud:06}. Then there exists a $\Delta$-net 
in $W_S$ in the $l_{\infty}$-metric such that
\begin{equation}
|\mathcal N|\le \left(\frac{C_{7.1}}{\Delta\sqrt{np_n}}\right)^n\eta^{c_{7.1}np_n}.
\end{equation}
\end{lem}
\begin{thm}\label{thm7.3}(Theorem 7.3 in \cite{rud:06})
There exists an absolute constant $Q_0$ with the following property. Let 
$\Delta>C_{7.3}(np_n)^{-\frac32}$, where $C_{7.3}=\max\{c_{4.2},\overline C_{7.1}\}$. Denote by 
$\Omega_{\Delta}$ the event that there exists a vector $\bold x\in V_S$ of $(\Delta,Q_0)$-singular profile
such that $\|\bold X^{\varepsilon}\bold x\|\le \frac{\Delta \sqrt{p_n}}{2}n$. Then
\begin{equation}
\Pr\{\Omega_{\Delta}\}\le 3\exp\{-np_n\}.
\end{equation}
\end{thm}
To prove Theorem \ref{thm1.1} we combine the probability estimates of the previous sections. Let $\gamma>\frac{c_{1.1}}{\sqrt{np_n}}$ where the constant $c_{1.1}$ will be chosen later. Define the exceptional sets:
\begin{align}
\Omega_0&=\{\omega\ |\ \|\bold X^{\varepsilon}\|>C_{2.3}\sqrt{np_n}\},\notag\\ 
\Omega_P&=\{\omega\ |\ \text{there exists }\bold x\in V_P \|\bold X^{\varepsilon}\bold x\|\le
C_{4.1}\sqrt{np_n}\}\notag.
\end{align}
Let $Q_0$ be the number defined in Theorem \ref{thm7.3}. Set 
\begin{equation}
\Delta=\frac{\gamma}{2C_{6.1}Q_0np_n}.
\end{equation}
The assumption on $\gamma$ implies $\Delta\ge C_{7.3}(np_n)^{-\frac32}$ if we set
$c_{1.1}=2C_{6.1}Q_0C_{7.3}$. Denote by $W_S$ the set of vectors of the $(\Delta,Q_0)$-singular profile and by $W_R$ the set of vectors of the $(\Delta,Q_0)$-regular profile. Set
\begin{align}
\Omega_S&=\left\{\omega\ \Big|\ \text{there exists } \bold x\in W_S \|\bold X^{\varepsilon}\bold x\|
\le \frac{\Delta\sqrt p_n}2n=\frac1{4C_{6.1Q_0}\sqrt{p_n}}\gamma\right\}\notag\\
\Omega_R&=\left\{\omega\ \Big|\ \text{there exists } \bold x\in W_R \|\bold X^{\varepsilon}\bold x\|
\le \frac{\Delta\sqrt{p_n}}{2\sqrt n}=\frac1{4C_{6.1Q_0}}\gamma n^{-\frac32}p_n^{-1}\right\}.
\end{align}
By Theorem \ref{thm7.3} $\Pr\{\Omega_S\}\le 3\exp\{-np_n\}$, and by Theorem \ref{thm6.2}
$\Pr\{\Omega_R\}\le Cn\Delta/p_n^{\frac12}\le C\gamma p_n^{-\frac32}$.
Choosing $\gamma =\frac c{\sqrt{np_n}}$,
we conclude the proof. 
\end{proof}


\end{document}